\documentclass{amsart}
\usepackage{amssymb,amsmath,color}
\usepackage{tikz}
\usetikzlibrary{shapes,arrows,backgrounds,fit,patterns}
\usetikzlibrary{positioning}
\usepackage{mathrsfs}

\tikzstyle{vertex} = [fill,shape=circle,node distance=80pt]
\tikzstyle{edge} = [fill,opacity=.5,fill opacity=.5,line cap=round, line join=round, line width=50pt]

\newtheorem{theorem}{Theorem}[section]

\newtheorem{lemma}[theorem]{Lemma}
\newtheorem{corollary}[theorem]{Corollary}
\newtheorem{proposition}[theorem]{Proposition}
\theoremstyle{definition}
\newtheorem{definition}[theorem]{Definition}

\theoremstyle{remark}

\date{\today}
\title{The Tur\'an Polytope}

\author{Annie Raymond}
\address{Department of Mathematics, University of Washington, Box
  354350, Seattle, WA 98195, USA} \email{raymonda@uw.edu}

\begin{document}

\begin{abstract}
The Tur\'an hypergraph problem asks to find the maximum number of $r$-edges in a $r$-uniform hypergraph on $n$ vertices that does not contain a clique of size $a$. When $r=2$, i.e., for graphs, the answer is well-known and can be found in Tur\'an's theorem. However, when $r\geq 3$, the problem remains open. We model the problem as an integer program and call the underlying polytope the Tur\'an polytope. We draw parallels between the latter and the stable set polytope: we show that generalized and transformed versions of the web and wheel inequalities are also facet-defining for the Tur\'an polytope. We also show clique inequalities and what we call doubling inequalities are facet-defining when $r=2$. These facets lead to a simple new polyhedral proof of Tur\'an's theorem.
\end{abstract}

\maketitle

\section{Introduction}

Mantel's theorem \cite{Mantel1907}, one of the earliest theorems in combinatorics, states that the maximum number of edges in a graph on $n$ vertices without any triangles is $\lfloor \frac{n^2}{4}\rfloor$, and that the maximum is attained only on complete bipartite graphs with parts of size $\lfloor \frac{n}{2} \rfloor$ and $\lceil \frac{n}{2} \rceil$.

Tur\'an \cite{Turan1941} later generalized this theorem by showing that the maximum number of edges in a graph on $n$ vertices without any clique of size $a$ is at most $(1-\frac{1}{a-1})\frac{n^2}{2}$, and that this maximum is attained solely on complete $(a-1)$-partite graphs with parts of size as equal as possible.

Since then, many different proofs of this theorem have been found using different techniques. We present a new polyhedral proof in Section~\ref{polyproof} by modeling the Tur\'an problem as an integer program. We prove in Theorem \ref{polyhedralproofturan} that the maximum number of edges in an $a$-clique free graph on $n$ vertices is exactly

 $$\left\lfloor \frac{n}{n-2} \left \lfloor \frac{n-1}{n-3} \left \lfloor \frac{n-2}{n-4} \left \lfloor \cdots \left \lfloor \frac{a+2}{a} \left \lfloor \frac{a+1}{a-1} \cdot \left(\binom{a}{2} -1\right) \right \rfloor \right \rfloor \cdots \right \rfloor \right \rfloor \right \rfloor \right \rfloor.$$

Tur\'an's theorem was later generalized in different ways. We turn ourselves to a generalization that Tur\'an himself introduced: the generalization to hypergraphs. The goal is to find the maximum number of $r$-hyperedges in a $r$-uniform hypergraph (i.e., a hypergraph for which every edge is composed of $r$ vertices) on $n$ vertices that does not contain any $r$-uniform hyperclique of size $a$ (i.e., a $r$-uniform hypergraph on $a$ vertices where every set of $r$ vertices forms an edge). 

This problem remains unsolved to this day. Even in the case when $a=4$ and $r=3$, the problem is still open. Tur\'an conjectured that in this case, the maximum number of edges is $(\frac{5}{9}+O(\frac{1}{n}))\binom{n}{3}$. The best known bound, 0.561666, is due to Razborov \cite{Razborov2010} who used flag algebra calculus. Moreover, Razborov also proved that if one additionally forbids graphs on four vertices that span exactly one $3$-edge as induced subgraphs, then the maximum hyperedge density is indeed $\frac{5}{9}$ as $n\rightarrow \infty$. 

Our integer program can easily be extended to the the Tur\'an hypergraph problem. We study some inequalities valid for the underlying polytope and show that they are facet-inducing under the right conditions. The inequalities we consider have a nice combinatorial flavor and draw a parallel between the Tur\'an polytope and the stable set polytope. Interestingly enough, we observe that the facets we study do not get dominated as $n$ increases, thus suggesting that the Tur\'an polytope is very complex. 

\subsection{Some Notation}

By $[n]$, we denote the set $\{1,2,\ldots,n\}$. In general, for a graph $G=(V,E)$ with $|V|=n$, we assume that $V=[n]$.  For a graph $G$, $V(G)$ and $E(G)$ represent respectively the vertex set and edge set of $G$. All of the graphs we consider are undirected, and we use $(v_1,\ldots, v_r)$ to denote an (undirected) edge formed by vertices $v_1$ through $v_r$. For $G=(V,E)$, $\delta(v)$ for some vertex $v\in V$ corresponds to the set of edges in $E$ that are adjacent to $v$ and $d(v):=|\delta(v)|$ is the \emph{degree} of $v$. Furthermore, for $S\subseteq V$, we let $E[S]$ be the set of edges induced by $S$ in $G$. We say a graph is \emph{$H$-free} or \emph{Tur\'an} if it does not contain an induced subgraph isomorphic to $H$. We let $K^r_n$ represent the complete $r$-uniform hypergraph on $n$ vertices, and we let $K_n:=K_n^2$.

\subsection{Previous Work}

The literature is simply too big to be included here. A whole paper could be written on the work done on Tur\'an-type problems for hypergraphs, and a whole book on Tur\'an-type problems in general. In fact, a whole paper \emph{has} been written on Tur\'an-type problems for hypergraphs: the excellent survey \cite{Keevash2011} by Peter Keevash. I am sitting on the edge of my seat waiting for someone to write a book on Tur\'an-type problems in general. In the meantime, one can consult \emph{Extremal Graph Theory} \cite{Bollobas1978} by Bollob\'as for results published before the mid-eighties. Together, these two bodies of work cover most of what is known. 

\section{Polyhedral Proof of Tur\'an's Theorem}\label{polyproof}

In this section, we give a polyhedral proof of Tur\'an's Theorem.

\begin{theorem}[Tur\'an, 1941]
The maximum number of edges in an $a$-clique-free graph on $n$ vertices is at most $(1-\frac{1}{a-1})\frac{n^2}{2}$. 
\end{theorem}

\subsection{Model}

\begin{definition}
Let $T(G,a,r)$ be the convex hull of the characteristic vectors of all edge sets $F\subseteq E(K^r_n)$ that contain no clique of size $a$.  Let $\mathscr{Q}^i_H$ be the set of all $r$-uniform hypercliques of size $i$, $Q^i$, in some $r$-uniform hypergraph $H$. Then

$$T(G,a,r)=\textup{conv}\left(\left\{x\in \{0,1\}^{\binom{n}{r}}| \sum_{e\in E(Q^a)} x_e \leq \binom{a}{r}-1 \quad \forall Q^a\in \mathscr{Q}^a_{G}\right\}\right).$$
We call $T(G,a,r)$ the \emph{Tur\'an polytope of $G$}. We denote $\max\{\sum_{e\in E(G)} x_e| x\in T(G,a,r)\}$ by $ex(G,a,r)$. We let $T(n,a,r):=T(K^r_n,a,r)$ and $ex(n,a,r):=ex(K^r_n,a,r)$. Furthermore, we let
\begin{displaymath}                                                                                                                                             \begin{array}{rlll}
Q(G,a,r)=\Big\{x\in \mathbb{R}^{\binom{n}{r}}|&\sum_{e\in E(Q^i)} x_e \leq ex(i,a,r) &\forall Q^i\in \mathscr{Q}^i_{G}, \forall a\leq i \leq n-1 & \\
& 0\leq x_e \leq 1 & \forall e\in E(G) & \Big\}                                                                                                                                  
\end{array}
\end{displaymath}
be the \emph{clique-relaxation of the Tur\'an polytope of $G$}. Again we let $Q(n,a,r):=Q(K^r_n,a,r)$. We call $\sum_{e\in E(Q^i)} x_e \leq ex(i,a,r)$ clique inequalities, $x_e \geq 0$ non-negativity inequalities, and $x_e \leq 1$ edge inequalities. 
\end{definition}

\subsection{Proof}

To prove Tur\'an's theorem, we need to show that $ex(n,a,2) \leq (1-\frac{1}{a-1})\frac{n^2}{2}$. To calculate $ex(n,a,2)$, we must consider two types of inequalities valid for $T(n,a,2)$.

\begin{definition}
Let $t_a^{i+1}= \lfloor \frac{i+1}{i-1} t_a^i \rfloor$ for $i\geq a$, and let $t_a^{a}=\binom{a}{2}-1$.
\end{definition}

\begin{lemma}
The inequality $\sum_{e\in E[S]} x_e \leq t_a^{|S|}$, where $E[S]$ consists of the edges of $K_n$ induced by the vertex set $S$, is valid on $T(n,a,2)$ for all $S\subseteq [n]$ such that $a\leq |S| \leq n$. Furthermore, so is $ex(n,a,2) \leq \left\lfloor \frac{n}{n-2} ex(n-1,a,2) \right\rfloor$.

\end{lemma}
\proof
We proceed by induction on the size of $|S|$. The base case $|S|=a$, corresponding to $\sum_{e\in E[S]} x_e \leq t_a^a=\binom{a}{2}-1$, is clear since the integer program for $T(n,a,2)$ contains these inequalities. Now assume that we have already shown the hypothesis for any $S$ such that $a\leq |S| \leq j$. Consider any set $S\in [n]$ such that $|S|=j+1$. Add up the inequalities corresponding to all $j+1$ sets $T$ of size $|S|-1$ contained in $S$. Notice that each edge in $E[S]$ is in $j-1$ of these inequalities, thus yielding that

$$\sum_{e\in E[S] }(j-1)x_e \leq (j+1) t_a^j,$$
is a valid inequality for $T(n,a,2)$ since it was produced as a conic combination of valid inequalities. Since we know that $\sum_{e\in E[S]}x_e$ is an integer for any $x\in T(n,a,2)$, 
$$\sum_{e\in E[S]}x_e \leq \left\lfloor\frac{j+1}{j-1} t_a^j\right\rfloor$$
is also a valid inequality for $T(n,a,2)$. Note that this argument also implies the upperbound \begin{equation}\label{eq1}ex(n,a,2) \leq \left\lfloor \frac{n}{n-2} ex(n-1,a,2) \right\rfloor.\end{equation}
\qed

Note that the Chv\'atal-Gomory cutting plane procedure applied to the linear relaxation $Q(n,a,2)$ would produce the inequality $\sum_{e\in E[S]} x_e \leq t_a^{|S|}$ at the latest in the $(|S|-a)$th round. Moreover, observe that, by letting $S=[n]$, we obtain that the inequality $\sum_{e\in E[S]} x_e=\sum_{e \in E(K_n)} x_e \leq t_a^n$ is valid for $T(n,a,2)$. By definition, this yields the upper bound $$ex(n,a,2) \leq t_a^n = \left\lfloor \frac{n}{n-2} \left \lfloor \frac{n-1}{n-3} \left \lfloor \frac{n-2}{n-4} \left \lfloor \cdots \left \lfloor \frac{a+2}{a} \left \lfloor \frac{a+1}{a-1} \cdot \left(\binom{a}{2} -1\right) \right \rfloor \right \rfloor \cdots \right \rfloor \right \rfloor \right \rfloor \right \rfloor.$$

To show that $ex(n,a,2)\geq t_a^n$ also, we show that another type of inequality is valid for $T(n,a,2)$. 

\begin{lemma}
The inequality $\sum_{e\in \delta(v)} 2x_e + \sum_{e\in E(K_n)\backslash \delta(v)} x_e \leq ex(n+1,a,2)$ for any vertex $v\in [n]$ is valid on $T(n,a,2)$. Furthermore, so is  $\sum_{e\in E(K_n)} x_e \leq \left\lfloor \frac{n}{n+2} ex(n+1,a,2)\right\rfloor$.

\end{lemma}

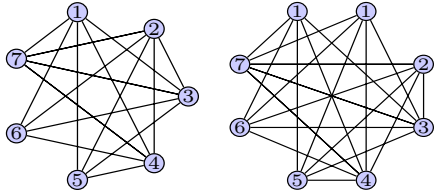
\begin{figure}[h]\label{doubling}
\begin{center}
\resizebox{3cm}{2.5cm}{
\begin{tikzpicture}
    % Define the heptagon coordinates
    \node [style={draw, fill=blue!20, circle, inner sep=0.3pt, minimum size=3pt}] (A) at (-0.76,1.54) {\tiny{7}};
    \node [style={draw, fill=blue!20, circle, inner sep=0.3pt, minimum size=3pt}] (B) at (-0.76,0.69) {\tiny{6}};
    \node [style={draw, fill=blue!20, circle, inner sep=0.3pt, minimum size=3pt}] (C) at (-0.10,0.16) {\tiny{5}};
    \node [style={draw, fill=blue!20, circle, inner sep=0.3pt, minimum size=3pt}] (D) at (0.73,0.35) {\tiny{4}};
    \node [style={draw, fill=blue!20, circle, inner sep=0.3pt, minimum size=3pt}] (E) at (1.1,1.11) {\tiny{3}};
    \node [style={draw, fill=blue!20, circle, inner sep=0.3pt, minimum size=3pt}] (F) at (0.73,1.88) {\tiny{2}};
    \node [style={draw, fill=blue!20, circle, inner sep=0.3pt, minimum size=3pt}] (G) at (-0.10,2.07) {\tiny{1}};

\draw (A)--(D)--(B)--(E)--(C)--(D)--(A)--(E)--(F)--(D)--(G)--(E)--(A)--(F)--(B)--(G)--(C)--(F)--(A)--(G);
\end{tikzpicture}
} \resizebox{3cm}{2.5cm}{\begin{tikzpicture}[scale=0.38]
    % Define the heptagon coordinates
    \node [style={draw, fill=blue!20, circle, inner sep=0.3pt, minimum size=3pt}] (A) at (-2.67,1) {\tiny{7}};
    \node [style={draw, fill=blue!20, circle, inner sep=0.3pt, minimum size=3pt}] (B) at (-2.67,-1) {\tiny{6}};
    \node [style={draw, fill=blue!20, circle, inner sep=0.3pt, minimum size=3pt}] (C) at (-1,-2.67) {\tiny{5}};
    \node [style={draw, fill=blue!20, circle, inner sep=0.3pt, minimum size=3pt}] (D) at (1,-2.67) {\tiny{4}};
    \node [style={draw, fill=blue!20, circle, inner sep=0.3pt, minimum size=3pt}] (E) at (2.67,-1) {\tiny{3}};
    \node [style={draw, fill=blue!20, circle, inner sep=0.3pt, minimum size=3pt}] (F) at (2.67,1) {\tiny{2}};
    \node [style={draw, fill=blue!20, circle, inner sep=0.3pt, minimum size=3pt}] (G) at (1,2.67) {\tiny{1}};
    \node [style={draw, fill=blue!20, circle, inner sep=0.3pt, minimum size=3pt}] (H) at (-1,2.67) {\tiny{1}};

\draw (A)--(D)--(B)--(E)--(C)--(D)--(A)--(E)--(F)--(D)--(G)--(E)--(A)--(F)--(B)--(G)--(C)--(F)--(A)--(G);
\draw (H)--(D);
\draw (H)--(E);
\draw (H)--(B);
\draw (H)--(C);
\draw (H)--(A);
\end{tikzpicture}
}
\caption{Copying vertex $1$ of a $4$-clique-free graph.}
\end{center}
\end{figure}

\proof
Take any $a$-clique-free graph $G$ with $n$ vertices. Fix a vertex $v$ and consider the graph $G'=([n+1], E')$ where $E'=E\cup \{(i,n+1)|(i,v)\in E\}$ (see Figure~\ref{doubling}). Note that $G'$ is also $a$-clique-free since any $a$ vertices that contains at most one of $n$ and $n+1$ cannot form a clique, otherwise $G$ would also have contained a clique, and no clique can contain both $n$ and $n+1$ since they don't form an edge.  The inequality $\sum_{e\in \delta(v)} 2x_e + \sum_{e\in E(K_n)\backslash \delta(v)} x_e \leq ex(n+1,a,2)$ is derived straight from that fact, and is thus valid on $T(n,a,2)$. 

Now, add up the $n$ inequalities corresponding to each $v\in [n]$. Note that each edge has coefficient of one in $n-2$ inequalities, and coefficient of two in two inequalities. We thus obtain that

$$\sum_{e\in E(K_n)} (n+2)x_e \leq n \cdot ex(n+1,a,2),$$
is a valid inequality for $T(n,a,2)$ since it was produced as a conic combination of valid inequalities. Since we know $\sum_{e\in E(K_n)}x_e$ is an integer for any $x\in T(n,a,2)$,
 $$\sum_{e\in E(K_n)} x_e \leq \left\lfloor \frac{n}{n+2} ex(n+1,a,2)\right\rfloor$$
is also a valid inequality for $T(n,a,2)$
\qed

Note that this yields the upper bound \begin{equation}\label{eq3}ex(n,a,2) \leq \left\lfloor \frac{n}{n+2} ex(n+1,a,2)\right\rfloor.\end{equation}

\begin{theorem}\label{polyhedralproofturan}
The equation $$ex(n,a,2)=\left\lfloor \frac{n}{n-2} \left \lfloor \frac{n-1}{n-3} \left \lfloor \frac{n-2}{n-4} \left \lfloor \cdots \left \lfloor \frac{a+2}{a} \left \lfloor \frac{a+1}{a-1} \cdot \left(\binom{a}{2} -1\right) \right \rfloor \right \rfloor \cdots \right \rfloor \right \rfloor \right \rfloor \right \rfloor$$ holds.
\end{theorem}

\proof
Putting inequalities (\ref{eq1}) and (\ref{eq3}) together, we get that \begin{equation}\label{together}ex(n+1,a,2) \leq \left\lfloor \frac{n+1}{n-1} \cdot ex(n,a,2) \right\rfloor \leq \left\lfloor \frac{n+1}{n-1} \left \lfloor \frac{n}{n+2} \cdot ex(n+1,a,2) \right \rfloor \right\rfloor.\end{equation} We now show that \begin{equation*}\left\lfloor \frac{n+1}{n-1} \left \lfloor \frac{n}{n+2} \cdot ex(n+1,a,2) \right \rfloor \right\rfloor \leq ex(n+1,a,2),\end{equation*} thus turning (\ref{together}) into an equation and proving our claim.

Suppose not. Then \begin{equation*}ex(n+1,a,2) + 1 \leq \left\lfloor \frac{n+1}{n-1} \left \lfloor \frac{n}{n+2} \cdot ex(n+1,a,2) \right \rfloor \right\rfloor,\end{equation*} which implies that: \begin{multline*}ex(n+1,a,2) + 1 \leq \left\lfloor \frac{n+1}{n-1} \cdot \frac{n}{n+2} \cdot ex(n+1,a,2)  \right\rfloor \\ = ex(n+1,a,2) + \left \lfloor \frac{2}{(n-1)\cdot(n+2)} \cdot ex(n+1,a,2) \right\rfloor.\end{multline*} Certainly, this means that\begin{align*} 1 &\leq \frac{2}{(n-1)\cdot(n+2)} \cdot ex(n+1,a,2) \\ &\leq \frac{2}{(n-1)\cdot(n+2)} \cdot \left\lfloor \frac{n+1}{n-1} \left \lfloor \frac{n}{n-2}  \left \lfloor \cdots \left \lfloor \frac{a+1}{a-1} \cdot \left(\binom{a}{2} -1\right) \right \rfloor \cdots \right \rfloor \right \rfloor \right \rfloor \\ &\leq \frac{2}{(n-1)\cdot(n+2)} \cdot \frac{n+1}{n-1} \cdot \frac{n}{n-2}  \cdots \frac{a+1}{a-1} \cdot \left(\binom{a}{2} -1\right) \\ &= \frac{(n+1) \cdot n \cdot (a+1) \cdot (a-2)}{(n-1)\cdot(n+2) \cdot a \cdot (a-1)}. \end{align*} So we have that \begin{equation*} \frac{(n-1)\cdot (n+2)}{(n+1)\cdot n} \leq \frac{(a+1)\cdot(a-2)}{a\cdot (a-1)} .\end{equation*}  By simplifying this inequality, we obtain that \begin{equation*} 1 - \frac{2}{n\cdot (n+1)} \leq 1 - \frac{2}{a\cdot (a-1)}, \end{equation*} which implies that \begin{equation*} n\cdot (n+1) \leq a\cdot(a-1),\end{equation*} which is impossible since $n > a \geq 3$. We have reached a contradiction, thus proving our claim.
\qed

\begin{corollary}
The integrality gap between $Q(n,a,2)$ and $T(n,a,2)$ is less than one in the $\mathbf{1}$-direction.
\end{corollary}

\proof
In the proof of the last theorem, we saw that

$$\left \lfloor \max_{x\in Q(n,a,2)} \mathbf{1}x\right \rfloor = ex(n,a,2) = \max_{x\in T(n,a,2)} \mathbf{1}x.$$
The result follows.
\qed

To the best of our knowledge, the exact formula

$$ex(n,a,2)=\left\lfloor \frac{n}{n-2} \left \lfloor \frac{n-1}{n-3} \left \lfloor \frac{n-2}{n-4} \left \lfloor \cdots \left \lfloor \frac{a+2}{a} \left \lfloor \frac{a+1}{a-1} \cdot \left(\binom{a}{2} -1\right) \right \rfloor \right \rfloor \cdots \right \rfloor \right \rfloor \right \rfloor \right \rfloor$$
was not known before. In the form $(1-\frac{1}{a-1})\frac{n^2}{2}$, Tur\'an's bound calculates the number of edges in an complete $a-1$-partite graph where all parts have equal size. It is of course only possible to do so if $n=0 \mod a-1$, and so the Tur\'an bound is equal to $ex(n,a,2)$ only if $n=0 \mod a-1$. For all other cases, $ex(n,a,2) < \left(1-\frac{1}{a-1}\right)\frac{n^2}{2}$.

\section{Facets of the Tur\'an Polytope}\label{polystudy}

In this section, we investigate some facet classes of $T(n,a,r)$. Some of these facets can be seen as analogs of famous facets of the stable set polytope. We first make a few general polytopal remarks that will help us generalize our results, then we show that the inequalities we described in the previous section are facet-defining for $T(n,a,2)$ under certain conditions, and we end with a proof that generalizations of web and wheel inequalities are facet-inducing for $T(n,a,r)$ under certain conditions. 

\subsection{General polytopal considerations} 

\begin{proposition}
Any facet we will add to the defining inequalities of $T(G,a,r)$ will be of the form $\alpha^Tx\leq b$ where $\alpha_e \geq 0$ for all $e\in E(G)$ and $b>0$. 
\end{proposition}

\proof
We simply need to show that $T(G,a,r)$ is down-monotone in $\mathbb{R}^n_+$. This is indeed the case since any $y$ such that $0\leq y \leq x$ for all $x\in T(G,a,r)$ will also be in $T(G,a,r)$ since $y$ will respect all the constraints of $T(G,a,r)$.
\qed

 We now show that some facets of $T(G,a,r)$ can be lifted to $T(G',a,r)$ where $G\subset G'$ because of the following theorem of Padberg.

\begin{theorem}[Padberg, 1973]\label{padberglifting}
  Let $S \subseteq \{0,1\}^n$ be monotone (i.e., $y \leq x \in S$ implies $y \in S$), $P_S := \textrm{conv}(S)$ be full-dimensional, $I \subseteq \{1,2,\ldots,n\}$, $P_S(I) := P_S \cap \{x\in \mathbb{R}^n | x_i = 0 \ \forall i \in I\}$ and $x^I \in \mathbb{R}^n$ denote a vector with $x_i^I = 0$ for all $i\in I$. Suppose $\sum_{j\not\in I} \alpha_jx_j \leq \alpha_o$ with $\alpha_0 > 0$ defines a facet of $P_s(I)$ and $i\in I$. Define

$$\alpha_i := \alpha_0 - \max\{\sum_{j\not\in I} \alpha_jx_j^I | e_i + x^I \in S\}.$$
Here $e_i$ is a unit vector in $\mathbb{R}^n$ with the $i$th component equal to one. Then

$$\alpha_ix_i + \sum_{j \not\in I} \alpha_jx_j \leq a_0$$
defines a facet of $P_S(I\backslash \{i\})$.
\end{theorem}

\begin{corollary}\label{turanliftinghyper01}
Let $G$ and $H$ be two $r$-uniform hypergraphs such that $H\subseteq G$ and such that $ex(H_e,a,r)=ex(H,a,r)+1$ for every $e\in E(G)\backslash E(H)$ where $H_e=(V(H),E(H)\cup e)$. If $\sum_{e\in E(H)}x_e \leq ex(H,a,r)$ is a facet of $T(H,a,r)$, then it is also a facet of $T(G,a,r)$.
\end{corollary}

\proof
We first note that $T(G,a,r)$ is a full-dimensional monotone polytope since $\mathbf{0}$ as well as every unit vector is in $T(G,a,r)$. Moreover, we know that if the characteristic vector of an edge set $S$ is in $T(G,a,r)$, then it is $a$-clique-free and if we take a subset of $S$, then it is also $a$-clique-free, and so its characteristic vector will also be in $T(G,a,r)$. Thus, $T(G,a,r)$ is also monotone, and so the setup is similar to the one of the previous theorem.

We let $I:=E(G)\backslash E(H)$, and we now show that $T(H,a,r)=T(G,a,r)\cap\{x\in \mathbb{R}^{|E(G)|}|x_e=0 \ \  \forall e \in I\}$. Indeed, $T(H,a,r)$ contains every point $P$ that is in $T(G,a,r)$ for which $x_e=0$ for all $e\in I$ since such a point must be a convex combination of vertices of $T(G,a,r)$ for which $x_e=0$ for all $e\in I$ (if the combination contained a positive coefficient for some vertex of $T(G,a,r)$ for which there exists $x_{e'}=1$ for some $e'\in I$, then $P$ would have $x_{e'}>0$, a contradiction). Since vertices of $T(G,a,r)$ for which $x_e=0$ are also vertices of $T(H,a,r)$, $P$ is also in $T(H,a,r)$. We can show similarly that every point in $T(H,a,r)$ is in $T(G,a,r)\cap\{x\in \mathbb{R}^{|E(G)|}|x_e=0 \ \ \forall e \in I\}$.

Suppose $\sum_{e\in E(H)} x_e \leq ex(H,a,r)$ is a facet of $T(H,a,r)$. Take any $e'\in I$, and let $c_{e'}:=ex(H,a,r) - \max\{\sum_{e\not\in I} x_e|x\in T(H,a,r) \textrm{ and } e'\cup x \textrm{ is } a-\textrm{clique-free}\}$. Since $ex(H_{e'},a,r)=ex(H,a,r)+1$, we know there exists an $a$-clique-free set of edges $S$ in $H_{e'}$ of size $ex(H,a,r)+1$ which must contain $e'$ (else we would have $ex(H,a,r)=ex(H,a,r)+1$!), and $S\backslash\{e'\}$ is still $a$-clique-free and completely in $H$, so $c_{e'}=0$ since $|S\backslash\{e'\}|=ex(H,a,r)$. Thus by the previous theorem, $\sum_{e\in E(H)} x_e \leq ex(H,a,r)$ is also a facet of $T(H_{e'},a,r)$.

Adding another edge $e'' \in E(G)\backslash (E(H)\cup\{e'\})$ will yield the same argument that $0=c_0 - \max\{\sum_{e\not\in I\backslash\{e'\}} x_e|x\in T(H_{e'},a,r) \textrm{ and } e''\cup x \textrm{ is } a-\textrm{clique-free}\}$ since $ex(H_{e''},a,r)=ex(H,a,r)+1$, and so there is an edge set of size $ex(H,a,r)$ in $T(H,a,r)$ (and thus also in $T(H_{e'},a,r)$) that fulfills the requirements.

We can therefore add all of the edges in $E(G)\backslash E(H)$, and get that $\sum_{e\in E(H)}x_e \leq ex(H,a,r)$ is a facet of $T(G,a,r).$ 
\qed

\begin{corollary}\label{lifting}
Let $G$ and $H$ be two $r$-uniform hypergraphs such that $H\subseteq G$. If $\sum_{e\in E(H)}c_ex_e \leq c_0$ is a facet of $T(H,a,r)$ such that for every $e'\in E(G)\backslash E(H)$ there exists $x^*$ such that $\sum_{e\in E(H)}c_ex^*_e = c_0$ for which $x^*\cup e'$ is $a$-clique-free, then it is also a facet of $T(G,a,r)$.
\end{corollary}

\proof
As in the previous corollary, we have that the setup is similar to that of the Padberg lifting property. Here again, let $I:=E(G)\backslash E(H)$ and $c_{e'}:=c_0 - \max\{\sum_{e\not\in I} x_e|x\in T(H,a,r) \textrm{ and } e'\cup x \textrm{ is } a-\textrm{clique-free}\}$ for $e' \in I$. Then, since there exists $x^*$ such that $x^*\cup e'$ is $a$-clique-free and and such that $\sum_{e\in E(H)}c_ex^*_e = c_0$, we have that $c_{e'}=0$ for any $e' \in I$. As in the previous corollary, we can add all of the edges in $I$ one after the other without encountering problems, and so we get that $\sum_{e\in E(H)} c_ex_e\leq c_0$ is a facet for $T(G,a,r)$. 
\qed

This is good news, and bad news at the same time. Many facets we find will still be facets on higher-dimensional examples; our work on smaller graphs will often carry on to larger graphs.  However, the fact that all those facets remain and don't get dominated by others when we add more edges to the graph means that higher-dimensional polytopes will have many, many, many facets, and thus makes it very unlikely that a complete description can be found for them.

\subsection{Some facets of $T(n,a,2)$}

\subsubsection{Clique facets}

The most trivial class of facets of the stable set polytope are the clique inequalities: $\sum_{v\in Q} x_v \leq 1$ for every clique $Q$. Such inequalities are facets of the stable set polytope for inclusionwise maximal cliques. Clique inequalities are also our most trivial class of facets. 

\begin{theorem}\label{turancliquefacet}
The clique inequality 

$$\sum_{e\in E(K_n)} x_e \leq ex(n,a,2)$$

is facet-defining for $T(n,a,2)$ if $n\neq 0 \mod (a-1)$.
\end{theorem}

\proof
From Tur\'an's theorem, we know that any $a$-clique-free graph with $ex(n,a,2)$ edges is an $(a-1)$-partite complete graph with parts of size as equal as possible. Suppose there are $p_1$ parts of size $\lfloor \frac{n}{a-1} \rfloor$ and $p_2$ parts of size $\lceil \frac{n}{a-1} \rceil$ such that $p_1+p_2 = a-1$ and $p_1\cdot\lfloor\frac{n}{a-1} \rfloor + p_2\cdot\lceil\frac{n}{a-1} \rceil = n$.

Let $\alpha x \leq \beta$ be satisfied by all $x$ in the Tur\'an polytope with $\sum_{e\in E(K_n)} x_e = ex(n,a,2)$; then $\alpha(S)=\beta$ for each Tur\'an edge set $S$ with $|S| = ex(n,a,2)$, i.e., the optimal $(a-1)$-partite complete graphs we've just described. 

Take two adjacent edges in $K_n$, without loss of generality, $(n,1)$ and $(1,2)$. We want to show that $\alpha_{(1,n)}=\alpha_{(1,2)}$. Since $n\neq 0 \mod (a-1)$, we know there exist optimal Tur\'an edge sets in the clique such that vertex $n$ is in a part of size $\lfloor \frac{n}{a-1} \rfloor $ and vertices $1$ and $2$ are together in a part of size $\lceil \frac{n}{a-1} \rceil$ such that $\lfloor \frac{n}{a-1} \rfloor < \lceil \frac{n}{a-1} \rceil$. For example, consider the optimal Tur\'an solution $S_1$ where vertices $n-\lfloor \frac{n}{a-1} \rfloor + 1$ through $n$ form one part, and vertices $1$ through $\lceil \frac{n}{a-1} \rceil$ form another part. Fix the other vertices into a partition $P$ that makes the whole solution optimal.

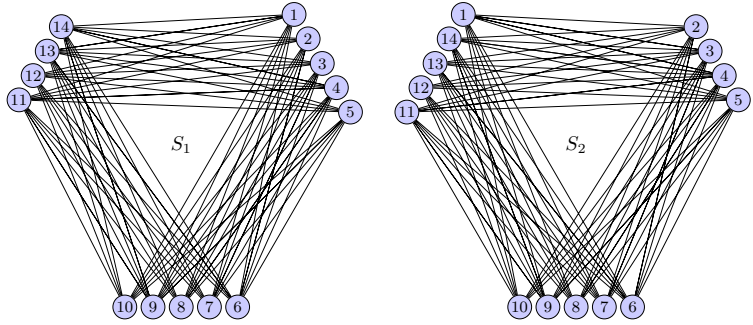
\begin{figure}[h]\label{figcliqueS1}
\begin{center}
\scalebox{.75}{\begin{tikzpicture}[scale=0.5]
    % Define the heptagon coordinates                                                                                                                                          
    \node [style={draw, fill=blue!20, circle, inner sep=0.3pt, minimum size=12pt}] (v1) at (4,4.619) {\footnotesize{1}};
    \node [style={draw, fill=blue!20, circle, inner sep=0.3pt, minimum size=12pt}] (v2) at (4.5,3.753) {\footnotesize{2}};
    \node [style={draw, fill=blue!20, circle, inner sep=0.3pt, minimum size=12pt}] (v3) at (5,2.88) {\footnotesize{3}};
    \node [style={draw, fill=blue!20, circle, inner sep=0.3pt, minimum size=12pt}] (v4) at (5.5,2.021) {\footnotesize{4}};
    \node [style={draw, fill=blue!20, circle, inner sep=0.3pt, minimum size=12pt}] (v5) at (6,1.155) {\footnotesize{5}};
    \node [style={draw, fill=blue!20, circle, inner sep=0.3pt, minimum size=12pt}] (v6) at (2,-5.76) {\footnotesize{6}};
    \node [style={draw, fill=blue!20, circle, inner sep=0.3pt, minimum size=12pt}] (v7) at (1,-5.76) {\footnotesize{7}};
    \node [style={draw, fill=blue!20, circle, inner sep=0.3pt, minimum size=12pt}] (v8) at (0,-5.76) {\footnotesize{8}};
    \node [style={draw, fill=blue!20, circle, inner sep=0.3pt, minimum size=12pt}] (v9) at (-1,-5.76) {\footnotesize{9}};
    \node [style={draw, fill=blue!20, circle, inner sep=0.3pt, minimum size=12pt}] (v10) at (-2,-5.76) {\footnotesize{10}};
    \node [style={draw, fill=blue!20, circle, inner sep=0.3pt, minimum size=12pt}] (v11) at (-5.75,1.588) {\footnotesize{11}};
    \node [style={draw, fill=blue!20, circle, inner sep=0.3pt, minimum size=12pt}] (v12) at (-5.25,2.45) {\footnotesize{12}};
    \node [style={draw, fill=blue!20, circle, inner sep=0.3pt, minimum size=12pt}] (v13) at (-4.75,3.317) {\footnotesize{13}};
    \node [style={draw, fill=blue!20, circle, inner sep=0.3pt, minimum size=12pt}] (v14) at (-4.25,4.186) {\footnotesize{14}};
    \node[draw=none,fill=none] at (0,0) {$S_1$};

\draw (v1) -- (v6) -- (v2) -- (v7)--(v3)--(v8)--(v4)--(v9)--(v5)--(v10)--(v1)--(v7)--(v3)--(v9)--(v5)--(v6)--(v2)--(v8)--(v4)--(v10)--(v2)--(v9)--(v1)--(v8)--(v5)--(v7)--(v4)--(v6)--(v3)--(v10);
\draw (v1)--(v11)--(v2)--(v12)--(v3)--(v13)--(v4)--(v14)--(v5)--(v13)--(v2)--(v14)--(v3)--(v11)--(v5)--(v12)--(v1)--(v13)--(v1)--(v14)--(v4)--(v11)--(v4)--(v12);
\draw (v6)--(v11)--(v7)--(v12)--(v8)--(v13)--(v9)--(v14)--(v10)--(v13)--(v7)--(v14)--(v8)--(v11)--(v10)--(v12)--(v6)--(v13)--(v6)--(v14)--(v9)--(v11)--(v9)--(v12);

    \node [style={draw, fill=blue!20, circle, inner sep=0.3pt, minimum size=12pt}] (w1) at (10,4.619) {\footnotesize{1}};
    \node [style={draw, fill=blue!20, circle, inner sep=0.3pt, minimum size=12pt}] (w14) at (9.5,3.753) {\footnotesize{14}};
    \node [style={draw, fill=blue!20, circle, inner sep=0.3pt, minimum size=12pt}] (w13) at (9,2.88) {\footnotesize{13}};
    \node [style={draw, fill=blue!20, circle, inner sep=0.3pt, minimum size=12pt}] (w12) at (8.5,2.021) {\footnotesize{12}};
    \node [style={draw, fill=blue!20, circle, inner sep=0.3pt, minimum size=12pt}] (w11) at (8,1.155) {\footnotesize{11}};
    \node [style={draw, fill=blue!20, circle, inner sep=0.3pt, minimum size=12pt}] (w6) at (16,-5.76) {\footnotesize{6}};
    \node [style={draw, fill=blue!20, circle, inner sep=0.3pt, minimum size=12pt}] (w7) at (15,-5.76) {\footnotesize{7}};
    \node [style={draw, fill=blue!20, circle, inner sep=0.3pt, minimum size=12pt}] (w8) at (14,-5.76) {\footnotesize{8}};
    \node [style={draw, fill=blue!20, circle, inner sep=0.3pt, minimum size=12pt}] (w9) at (13,-5.76) {\footnotesize{9}};
    \node [style={draw, fill=blue!20, circle, inner sep=0.3pt, minimum size=12pt}] (w10) at (12,-5.76) {\footnotesize{10}};
    \node [style={draw, fill=blue!20, circle, inner sep=0.3pt, minimum size=12pt}] (w5) at (19.75,1.588) {\footnotesize{5}};
    \node [style={draw, fill=blue!20, circle, inner sep=0.3pt, minimum size=12pt}] (w4) at (19.25,2.45) {\footnotesize{4}};
    \node [style={draw, fill=blue!20, circle, inner sep=0.3pt, minimum size=12pt}] (w3) at (18.75,3.317) {\footnotesize{3}};
    \node [style={draw, fill=blue!20, circle, inner sep=0.3pt, minimum size=12pt}] (w2) at (18.25,4.186) {\footnotesize{2}};
    \node[draw=none,fill=none] at (14,0) {$S_2$};

\draw (w1) -- (w6) -- (w2) -- (w7)--(w3)--(w8)--(w4)--(w9)--(w5)--(w10)--(w1)--(w7)--(w3)--(w9)--(w5)--(w6)--(w2)--(w8)--(w4)--(w10)--(w2)--(w9)--(w1)--(w8)--(w5)--(w7)--(w4)--(w6)--(w3)--(w10);
\draw (w11)--(w2)--(w12)--(w3)--(w13)--(w4)--(w14)--(w5)--(w13)--(w2)--(w14)--(w3)--(w11)--(w5)--(w12);
\draw (w14)--(w4)--(w11)--(w4)--(w12);
\draw (w6)--(w11)--(w7)--(w12)--(w8)--(w13)--(w9)--(w14)--(w10)--(w13)--(w7)--(w14)--(w8)--(w11)--(w10)--(w12)--(w6)--(w13)--(w6)--(w14)--(w9)--(w11)--(w9)--(w12);
\draw (w2)--(w1)--(w3)--(w1)--(w4)--(w1)--(w5);

\end{tikzpicture}}
\caption{Example of $S_1$ and $S_2$ for $n=14$ and $a=4$.}
\end{center}
\end{figure}

Now consider another optimal solution $S_2$ such that vertices $\lceil \frac{n}{a-1} \rceil + 1$ through $n-\lfloor \frac{n}{a-1} \rfloor$ are partitioned into $P$, vertices $n-\lfloor \frac{n}{a-1} \rfloor + 1$ through $n$ with vertex $1$ as well form a part of size $\lceil \frac{n}{a-1} \rceil$ and vertices $2$ through $\lceil\frac{n}{a-1}\rceil$ form another part of size $\lfloor \frac{n}{a-1} \rfloor$. This is thus the previous solution but with vertex $1$ moved to the other defined part. Note that we couldn't do that if $n=0 \mod (a-1)$ and still have an optimal solution after the move.

Since $\alpha(S_1)=\beta=\alpha(S_2)$, this implies that 

\begin{multline}\label{eqS12clique}
\alpha_{(1,n-\lfloor \frac{n}{a-1} \rfloor + 1)} + \alpha_{(1,n-\lfloor \frac{n}{a-1} \rfloor + 2)}+ \ldots + \alpha_{(1,n)} \\
= \alpha_{(1,2)} + \alpha_{(1,3)} + \ldots + \alpha_{(1,\lceil\frac{n}{a-1}\rceil)}.
\end{multline}

Now consider yet another optimal Tur\'an edge set $S_3$ such that again vertices $\lceil \frac{n}{a-1} \rceil + 1$ through $n-\lfloor \frac{n}{a-1} \rfloor$ are partitioned into $P$, vertices $n-\lfloor \frac{n}{a-1} \rfloor + 1$ through $n-1$ with vertex $2$ as well form a part of size $\lfloor \frac{n}{a-1} \rfloor$ and vertices $3$ through $\lceil\frac{n}{a-1}\rceil$ with vertices $1$ and $n$ form another part of size $\lceil \frac{n}{a-1} \rceil$. 

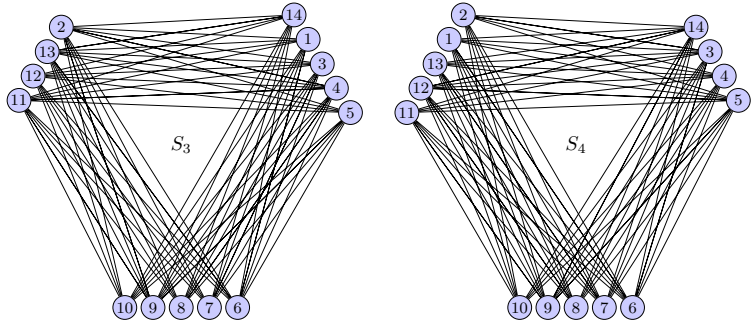
\begin{figure}[h]\label{figcliqueS3S4}
\begin{center}
\scalebox{.75}{\begin{tikzpicture}[scale=0.5]
    % Define the heptagon coordinates                                                                                                                                          
    \node [style={draw, fill=blue!20, circle, inner sep=0.3pt, minimum size=12pt}] (v1) at (4,4.619) {\footnotesize{14}};
    \node [style={draw, fill=blue!20, circle, inner sep=0.3pt, minimum size=12pt}] (v2) at (4.5,3.753) {\footnotesize{1}};
    \node [style={draw, fill=blue!20, circle, inner sep=0.3pt, minimum size=12pt}] (v3) at (5,2.88) {\footnotesize{3}};
    \node [style={draw, fill=blue!20, circle, inner sep=0.3pt, minimum size=12pt}] (v4) at (5.5,2.021) {\footnotesize{4}};
    \node [style={draw, fill=blue!20, circle, inner sep=0.3pt, minimum size=12pt}] (v5) at (6,1.155) {\footnotesize{5}};
    \node [style={draw, fill=blue!20, circle, inner sep=0.3pt, minimum size=12pt}] (v6) at (2,-5.76) {\footnotesize{6}};
    \node [style={draw, fill=blue!20, circle, inner sep=0.3pt, minimum size=12pt}] (v7) at (1,-5.76) {\footnotesize{7}};
    \node [style={draw, fill=blue!20, circle, inner sep=0.3pt, minimum size=12pt}] (v8) at (0,-5.76) {\footnotesize{8}};
    \node [style={draw, fill=blue!20, circle, inner sep=0.3pt, minimum size=12pt}] (v9) at (-1,-5.76) {\footnotesize{9}};
    \node [style={draw, fill=blue!20, circle, inner sep=0.3pt, minimum size=12pt}] (v10) at (-2,-5.76) {\footnotesize{10}};
    \node [style={draw, fill=blue!20, circle, inner sep=0.3pt, minimum size=12pt}] (v11) at (-5.75,1.588) {\footnotesize{11}};
    \node [style={draw, fill=blue!20, circle, inner sep=0.3pt, minimum size=12pt}] (v12) at (-5.25,2.45) {\footnotesize{12}};
    \node [style={draw, fill=blue!20, circle, inner sep=0.3pt, minimum size=12pt}] (v13) at (-4.75,3.317) {\footnotesize{13}};
    \node [style={draw, fill=blue!20, circle, inner sep=0.3pt, minimum size=12pt}] (v14) at (-4.25,4.186) {\footnotesize{2}};
    \node[draw=none,fill=none] at (0,0) {$S_3$}; 

\draw (v1) -- (v6) -- (v2) -- (v7)--(v3)--(v8)--(v4)--(v9)--(v5)--(v10)--(v1)--(v7)--(v3)--(v9)--(v5)--(v6)--(v2)--(v8)--(v4)--(v10)--(v2)--(v9)--(v1)--(v8)--(v5)--(v7)--(v4)--(v6)--(v3)--(v10);
\draw (v1)--(v11)--(v2)--(v12)--(v3)--(v13)--(v4)--(v14)--(v5)--(v13)--(v2)--(v14)--(v3)--(v11)--(v5)--(v12)--(v1)--(v13)--(v1)--(v14)--(v4)--(v11)--(v4)--(v12);
\draw (v6)--(v11)--(v7)--(v12)--(v8)--(v13)--(v9)--(v14)--(v10)--(v13)--(v7)--(v14)--(v8)--(v11)--(v10)--(v12)--(v6)--(v13)--(v6)--(v14)--(v9)--(v11)--(v9)--(v12);

    \node [style={draw, fill=blue!20, circle, inner sep=0.3pt, minimum size=12pt}] (w1) at (10,4.619) {\footnotesize{2}};
    \node [style={draw, fill=blue!20, circle, inner sep=0.3pt, minimum size=12pt}] (w2) at (9.5,3.753) {\footnotesize{1}};
    \node [style={draw, fill=blue!20, circle, inner sep=0.3pt, minimum size=12pt}] (w3) at (9,2.88) {\footnotesize{13}};
    \node [style={draw, fill=blue!20, circle, inner sep=0.3pt, minimum size=12pt}] (w4) at (8.5,2.021) {\footnotesize{12}};
    \node [style={draw, fill=blue!20, circle, inner sep=0.3pt, minimum size=12pt}] (w5) at (8,1.155) {\footnotesize{11}};
    \node [style={draw, fill=blue!20, circle, inner sep=0.3pt, minimum size=12pt}] (w6) at (16,-5.76) {\footnotesize{6}};
    \node [style={draw, fill=blue!20, circle, inner sep=0.3pt, minimum size=12pt}] (w7) at (15,-5.76) {\footnotesize{7}};
    \node [style={draw, fill=blue!20, circle, inner sep=0.3pt, minimum size=12pt}] (w8) at (14,-5.76) {\footnotesize{8}};
    \node [style={draw, fill=blue!20, circle, inner sep=0.3pt, minimum size=12pt}] (w9) at (13,-5.76) {\footnotesize{9}};
    \node [style={draw, fill=blue!20, circle, inner sep=0.3pt, minimum size=12pt}] (w10) at (12,-5.76) {\footnotesize{10}};
    \node [style={draw, fill=blue!20, circle, inner sep=0.3pt, minimum size=12pt}] (w11) at (19.75,1.588) {\footnotesize{5}};
    \node [style={draw, fill=blue!20, circle, inner sep=0.3pt, minimum size=12pt}] (w12) at (19.25,2.45) {\footnotesize{4}};
    \node [style={draw, fill=blue!20, circle, inner sep=0.3pt, minimum size=12pt}] (w13) at (18.75,3.317) {\footnotesize{3}};
    \node [style={draw, fill=blue!20, circle, inner sep=0.3pt, minimum size=12pt}] (w14) at (18.25,4.186) {\footnotesize{14}};
    \node[draw=none,fill=none] at (14,0) {$S_4$}; 

\draw (w1) -- (w6) -- (w2) -- (w7)--(w3)--(w8)--(w4)--(w9)--(w5)--(w10)--(w1)--(w7)--(w3)--(w9)--(w5)--(w6)--(w2)--(w8)--(w4)--(w10)--(w2)--(w9)--(w1)--(w8)--(w5)--(w7)--(w4)--(w6)--(w3)--(w10);
\draw (w1)--(w11)--(w2)--(w12)--(w3)--(w13)--(w4)--(w14)--(w5)--(w13)--(w2)--(w14)--(w3)--(w11)--(w5)--(w12)--(w1)--(w13)--(w1)--(w14)--(w4)--(w11)--(w4)--(w12);
\draw (w6)--(w11)--(w7)--(w12)--(w8)--(w13)--(w9)--(w14)--(w10)--(w13)--(w7)--(w14)--(w8)--(w11)--(w10)--(w12)--(w6)--(w13)--(w6)--(w14)--(w9)--(w11)--(w9)--(w12);
\end{tikzpicture}}
\caption{Example of $S_3$ and $S_4$ for $n=14$ and $a=4$.}
\end{center}
\end{figure}

Finally, we consider one last optimal Tur\'an solution $S_4$, again with vertices $\lceil \frac{n}{a-1} \rceil + 1$ through $n-\lfloor \frac{n}{a-1} \rfloor$ are partitioned into $P$, vertices $n-\lfloor \frac{n}{a-1} \rfloor + 1$ through $n-1$ with vertices $1$ and $2$ as well form a part of size $\lceil \frac{n}{a-1} \rceil$ and vertices $3$ through $\lceil \frac{n}{a-1}\rceil$ with vertex $n$ form another part of size $\lfloor \frac{n}{a-1} \rfloor$ (thus the previous solution but with vertex $1$ moved to the other defined part). 

Again, since $\alpha(S_3)=\beta=\alpha(S_4)$, this implies that 

\begin{multline}\label{eqS34clique}
\alpha_{(1,2)}+\alpha_{(1,n-\lfloor \frac{n}{a-1} \rfloor + 1)} + \alpha_{(1,n-\lfloor \frac{n}{a-1} \rfloor + 2)}+ \ldots + \alpha_{(1,n-1)}\\
=\alpha_{(1,3)}  + \alpha_{(1,4)} + \ldots + \alpha_{(1,\lceil\frac{n}{a-1}\rceil)} + \alpha_{(1,n)}.
\end{multline}

By subtracting equation \ref{eqS34clique} from equation \ref{eqS12clique}, we get that

$$\alpha_{(1,n)} - \alpha_{(1,2)} = \alpha_{(1,2)} - \alpha_{(1,n)},$$

which implies that $\alpha_{(1,n)}=\alpha_{(1,2)}$. Since we chose those edges without loss of generality, we know that $\alpha_{e_1}=\alpha_{e_2}$ for any two adjacent edges $e_1$ and $e_2$. By applying this observation to all pairs of adjacent edges in the clique, we obtain that $\alpha_{e_3}=\alpha_{e_4}$ for any two edges $e_3, e_4$ of $K_n$. Note that it is clear that $\alpha_e >0$. Thus $\alpha$ is a positive scalar multiple of the left-hand side of the clique inequality, which is thus facet-defining when $n \neq 0 \mod (a-1)$.
\qed

\begin{corollary}\label{turancliquelifting}
If $G$ contains a clique of size $i$ with $i \neq 0 \mod (a-1)$, say $Q^i$, then the corresponding clique inequality is a facet of $T(G,a,2)$. 
\end{corollary}

\proof
Note that $ex(Q^i+e, a, 2)=ex(i,a,2)+1$ for all $e\in E(G)\backslash E(Q^i)$ since any edge $e$ in $G$ but not in the clique can be added to any optimal edge set of the clique without forming an $a$-clique since other edges containing some vertices of both $e$ and $Q^i$ are missing. Thus, by Corollary \ref{lifting}, these clique inequalities are facets of $T(G,a,2)$. In particular, they are facets of $T(n,a,2)$ for $n\geq i$. 
\qed

\subsubsection{Doubling Facets}
In the polyhedral proof of the Tur\'an theorem, we introduced another type of valid inequality, which we call the \emph{doubling inequality}. We showed there that

$$\sum_{e\in \delta(v)} 2x_e + \sum_{e \in E(K_n)\backslash \delta(v)} x_e \leq ex(n+1,a,2)$$
is a valid inequality for $T(n,a,2)$ for any $v\in [n]$ since copying any vertex in a Tur\'an edge set gives an edge set that is also Tur\'an. We now show that this inequality is sometimes facet-defining.

\begin{theorem}\label{turandoublingfacet}
The doubling inequality 

$$\sum_{e\in \delta(v)} 2x_e + \sum_{e \in E(K_n)\backslash \delta(v)} x_e \leq ex(n+1,a,2)$$
is facet-defining for any $v\in [n]$ for $T(n,a,2)$ when $n=0 \mod (a-1)$ and with $n\geq 3(a-1)$.
\end{theorem}

\proof
We first observe that there are two types of edge sets that are tight with the doubling inequality for $n=0 \mod (a-1)$. The first type, which we call type I, is simply the optimal Tur\'an edge set of a clique, that is, a complete  $(a-1)$-partite graph with each part containing $\frac{n}{a-1}$ vertices. There are $\binom{\frac{n}{a-1}}{2}\cdot(\frac{n}{a-1})^2$ edges in such a solution, including $n-\frac{n}{a-1}$ edges that get doubled in the doubling inequality, thus the left-hand side of the inequality yield $\binom{\frac{n}{a-1}}{2}\cdot(\frac{n}{a-1})^2 + n-\frac{n}{a-1}$ which is equal to $ex(n+1,a,2)$, as desired.

The second type of optimal Tur\'an edge set, called type II, is given by the same construction, but with one vertex being moved to another part and the vertex $v$ that gets doubled being in the part that lost a vertex. That is, an $(a-1)$-partite graph with all parts containing $\frac{n}{a-1}$ vertices except for two parts, one containing $\frac{n}{a-1}-1$ vertices (including $v$) and one containing $\frac{n}{a-1}+1$ vertices. 

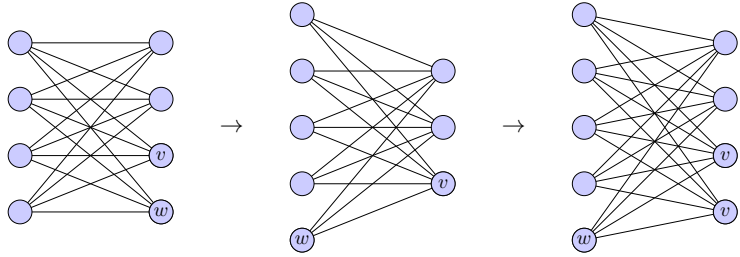
\begin{figure}[h]\label{figdoublingtypeii}
\begin{center}
\scalebox{0.75}{\begin{tikzpicture}[scale=1]
\foreach \x in {0.5,1.5,2.5,3.5}
\foreach \y in {0.5,1.5,2.5,3.5}
{\draw (0,\y) -- (2.5,\x);}
\foreach \x in {0.5,1.5,2.5,3.5}{
\node [style={draw, fill=blue!20, circle, inner sep=0.3pt, minimum size=12pt}] (l\x) at (0,\x) {};
\node [style={draw, fill=blue!20, circle, inner sep=0.3pt, minimum size=12pt}] (r\x) at (2.5,\x) {};
}

\node [style={draw, fill=blue!20, circle, inner sep=0.3pt, minimum size=12pt}] (v) at (2.5,0.5) {$w$};
\node [style={draw, fill=blue!20, circle, inner sep=0.3pt, minimum size=12pt}] (w) at (2.5,1.5) {$v$};

  \node[draw=none,fill=none] at (3.75,2) {\Large{$\mathbf{\rightarrow}$}};

\foreach \x in {0,1,2,3,4}
\foreach \y in {1,2,3}
{\draw (7.5,\y) -- (5,\x);}
\foreach \x in {0,1,2,3,4}{
\node [style={draw, fill=blue!20, circle, inner sep=0.3pt, minimum size=12pt}] (a\x) at (5,\x) {};}
\foreach \y in {1,2,3}{
\node [style={draw, fill=blue!20, circle, inner sep=0.3pt, minimum size=12pt}] (b\y) at (7.5,\y) {};}

\node [style={draw, fill=blue!20, circle, inner sep=0.3pt, minimum size=12pt}] (v) at (7.5,1) {$v$};
\node [style={draw, fill=blue!20, circle, inner sep=0.3pt, minimum size=12pt}] (w) at (5,0) {$w$};

  \node[draw=none,fill=none] at (8.75,2) {\Large{$\mathbf{\rightarrow}$}};
  
\foreach \x in {0,1,2,3,4}
\foreach \y in {0.5,1.5,2.5,3.5}
{\draw (12.5,\y) -- (10,\x);}
\foreach \x in {0,1,2,3,4}{
\node [style={draw, fill=blue!20, circle, inner sep=0.3pt, minimum size=12pt}] (a\x) at (10,\x) {};}
\foreach \y in {0.5,1.5,2.5,3.5}{
\node [style={draw, fill=blue!20, circle, inner sep=0.3pt, minimum size=12pt}] (b\y) at (12.5,\y) {};}

\node [style={draw, fill=blue!20, circle, inner sep=0.3pt, minimum size=12pt}] (v) at (12.5,1.5) {$v$};
\node [style={draw, fill=blue!20, circle, inner sep=0.3pt, minimum size=12pt}] (v) at (12.5,0.5) {$v$};
\node [style={draw, fill=blue!20, circle, inner sep=0.3pt, minimum size=12pt}] (w) at (10,0) {$w$};

\end{tikzpicture}}
\caption{Constructing an optimal solution of type II for the doubling inequality for $T(8,3,2)$}
\end{center}
\end{figure}

One can check that the number of edges in a type II solution is the same as in a type I solution by calculating how many edges we lose and gain, 
%. So by moving a vertex from one part to another in the optimal Tur\'an edge set for the clique to get the basis of a type II solution, the only edges that change are those adjacent to the vertex $w$ we move. In the original clique solution, the vertex $w$ we move was adjacent to $n-\frac{n}{a-1}$ edges, none of them adjacent to $v$, so forming a total of $n-\frac{n}{a-1}$ in the left-hand side of the doubling inequality. In the second type of solution, $w$ is adjacent to only $n-\frac{n}{a-1}-1$ edges, but one of them is adjacent to $v$, and so they form a total of $n-\frac{n}{a-1}$ in the left-hand side of the doubling inequality. So the total change is zero,
and that the second type of solution is also tight with the doubling inequality if $n=0\mod a-1$. Note that the type II construction is not always optimal for the doubling inequality when $n\neq 0 \mod a-1$, and observe that the type I constructions are already in the clique facet in these cases. 

Let $\alpha x \leq \beta$ be a facet of $T(n,a,2)$  satisfied by all $x$ in the Tur\'an polytope with $\sum_{e\in \delta(v)} 2x_e + \sum_{e\in E(K_n)\backslash \delta(v)} x_e = ex(n+1,a,2)$; then $\alpha(S)=\beta$ for both types of Tur\'an edge sets that we have just described.

Without loss of generality, let $v=\frac{n}{a-1}$. Take two adjacent edges that do not contain $v$, without loss of generality, $(n,1)$ and $(1,2)$. First consider the optimal solution $S_1$ of type I where vertices $i\cdot\frac{n}{a-1} + 1$ through $(i+1)\cdot\frac{n}{a-1}$ form a part for $0\leq i \leq a-2$. 

\begin{figure}[h]\label{figdoublingS1S2}
\begin{center}
\scalebox{.75}{\begin{tikzpicture}[scale=1]
\foreach \x in {0.5,1.5,2.5,3.5}
\foreach \y in {0.5,1.5,2.5,3.5}
{\draw (0,\y) -- (2.5,\x);}

\node [style={draw, fill=blue!20, circle, inner sep=0.3pt, minimum size=12pt}] (1) at (0,0.5) {1};
\node [style={draw, fill=blue!20, circle, inner sep=0.3pt, minimum size=12pt}] (2) at (0,1.5) {2};
\node [style={draw, fill=blue!20, circle, inner sep=0.3pt, minimum size=12pt}] (3) at (0,2.5) {3};
\node [style={draw, fill=blue!20, circle, inner sep=0.3pt, minimum size=12pt}] (4) at (0,3.5) {4};
\node [style={draw, fill=blue!20, circle, inner sep=0.3pt, minimum size=12pt}] (5) at (2.5,0.5) {5};
\node [style={draw, fill=blue!20, circle, inner sep=0.3pt, minimum size=12pt}] (6) at (2.5,1.5) {6};
\node [style={draw, fill=blue!20, circle, inner sep=0.3pt, minimum size=12pt}] (7) at (2.5,2.5) {7};
\node [style={draw, fill=blue!20, circle, inner sep=0.3pt, minimum size=12pt}] (8) at (2.5,3.5) {8};

  \node[draw=none,fill=none] at (1.25,0) {\Large{$S_1$}};

\foreach \x in {0,1,2,3,4}
\foreach \y in {1,2,3}
{\draw (5,\y) -- (7.5,\x);}

\node [style={draw, fill=blue!20, circle, inner sep=0.3pt, minimum size=12pt}] (v) at (7.5,0) {$1$};
\node [style={draw, fill=blue!20, circle, inner sep=0.3pt, minimum size=12pt}] (v) at (7.5,1) {$5$};
\node [style={draw, fill=blue!20, circle, inner sep=0.3pt, minimum size=12pt}] (v) at (7.5,2) {$6$};
\node [style={draw, fill=blue!20, circle, inner sep=0.3pt, minimum size=12pt}] (v) at (7.5,3) {$7$};
\node [style={draw, fill=blue!20, circle, inner sep=0.3pt, minimum size=12pt}] (v) at (7.5,4) {$8$};
\node [style={draw, fill=blue!20, circle, inner sep=0.3pt, minimum size=12pt}] (w) at (5,1) {$2$};
\node [style={draw, fill=blue!20, circle, inner sep=0.3pt, minimum size=12pt}] (w) at (5,2) {$3$};
\node [style={draw, fill=blue!20, circle, inner sep=0.3pt, minimum size=12pt}] (w) at (5,3) {$4$};

  \node[draw=none,fill=none] at (6.25,0) {\Large{$S_2$}};

\end{tikzpicture}}
\caption{$S_1$ and $S_2$ for $T(8,3,2)$}
\end{center}
\end{figure}
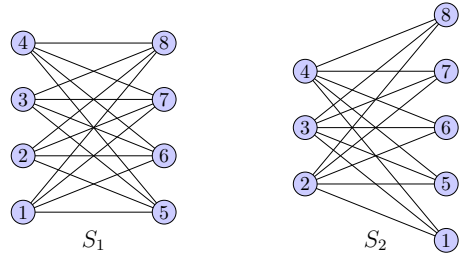

Now consider the solution of type II $S_2$ which is the same as $S_1$ but with vertex $1$ moved to the part containing vertices $n - \frac{n}{a-1} + 1$ through $n$. Since $\alpha(S_1)=\beta= \alpha(S_2)$, we have that

$$\alpha_{(1,n-\frac{n}{a-1} + 1)} + \alpha_{(1,n-\frac{n}{a-1} + 2)} + \ldots + \alpha_{(1,n)}= \alpha_{(1,2)} + \alpha_{(1,3)} + \ldots + \alpha_{(1,\frac{n}{a-1})}.$$

Consider now two new optimal solutions. First, we consider one of type I, say $S_3$, where vertices $i\cdot\frac{n}{a-1} + 1$ through $(i+1)\cdot\frac{n}{a-1}$ form a part for $1\leq i \leq a-3$, vertices $3$ through $\frac{n}{a-1}$ with vertices $n$ and $1$ form another part, and vertices $n-\frac{n}{a-1}+1$ through $n-1$ with vertex $2$ form a final part.

\begin{figure}[h]\label{figdoublingS3S4}
\begin{center}
\scalebox{.75}{\begin{tikzpicture}[scale=1]
\foreach \x in {0.5,1.5,2.5,3.5}
\foreach \y in {0.5,1.5,2.5,3.5}
{\draw (0,\y) -- (2.5,\x);}

\node [style={draw, fill=blue!20, circle, inner sep=0.3pt, minimum size=12pt}] (1) at (0,0.5) {8};
\node [style={draw, fill=blue!20, circle, inner sep=0.3pt, minimum size=12pt}] (2) at (0,1.5) {1};
\node [style={draw, fill=blue!20, circle, inner sep=0.3pt, minimum size=12pt}] (3) at (0,2.5) {3};
\node [style={draw, fill=blue!20, circle, inner sep=0.3pt, minimum size=12pt}] (4) at (0,3.5) {4};
\node [style={draw, fill=blue!20, circle, inner sep=0.3pt, minimum size=12pt}] (5) at (2.5,0.5) {2};
\node [style={draw, fill=blue!20, circle, inner sep=0.3pt, minimum size=12pt}] (6) at (2.5,1.5) {5};
\node [style={draw, fill=blue!20, circle, inner sep=0.3pt, minimum size=12pt}] (7) at (2.5,2.5) {6};
\node [style={draw, fill=blue!20, circle, inner sep=0.3pt, minimum size=12pt}] (8) at (2.5,3.5) {7};

  \node[draw=none,fill=none] at (1.25,0) {\Large{$S_3$}};

\foreach \x in {0,1,2,3,4}
\foreach \y in {1,2,3}
{\draw (5,\y) -- (7.5,\x);}

\node [style={draw, fill=blue!20, circle, inner sep=0.3pt, minimum size=12pt}] (v) at (7.5,0) {$1$};
\node [style={draw, fill=blue!20, circle, inner sep=0.3pt, minimum size=12pt}] (v) at (7.5,1) {$2$};
\node [style={draw, fill=blue!20, circle, inner sep=0.3pt, minimum size=12pt}] (v) at (7.5,2) {$5$};
\node [style={draw, fill=blue!20, circle, inner sep=0.3pt, minimum size=12pt}] (v) at (7.5,3) {$6$};
\node [style={draw, fill=blue!20, circle, inner sep=0.3pt, minimum size=12pt}] (v) at (7.5,4) {$7$};
\node [style={draw, fill=blue!20, circle, inner sep=0.3pt, minimum size=12pt}] (w) at (5,1) {$8$};
\node [style={draw, fill=blue!20, circle, inner sep=0.3pt, minimum size=12pt}] (w) at (5,2) {$3$};
\node [style={draw, fill=blue!20, circle, inner sep=0.3pt, minimum size=12pt}] (w) at (5,3) {$4$};

  \node[draw=none,fill=none] at (6.25,0) {\Large{$S_4$}};

\end{tikzpicture}}
\caption{$S_3$ and $S_4$ for $T(8,3,2)$}
\end{center}
\end{figure}
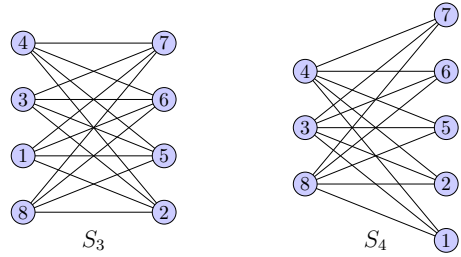

Second, we consider a final solution $S_4$ of type II which is the same as $S_3$ but with vertex $1$ moved to the part containing  vertices $n-\frac{n}{a-1}+1$ through $n-1$ as well as vertex $2$. Since $\alpha(S_3)=\beta=\alpha(S_4)$, we obtain that

\begin{multline*}
\alpha_{(1,2)}+\alpha_{(1,n-\frac{n}{a-1}+1)} + \alpha_{(1,n-\frac{n}{a-1}+2)}+\ldots \alpha_{(1,n-1)}\\
= \alpha_{(1,3)}+\alpha_{(1,4)} + \ldots + \alpha_{(1,\frac{n}{a-1})}+\alpha_{(1,n)}.
\end{multline*}
Subtracting these two equations, we get that

$$\alpha_{(1,n)} - \alpha_{(1,2)} = \alpha_{(1,2)}-\alpha_{(1,n)},$$
meaning that $\alpha_{(1,n)}=\alpha_{(1,2)}$. Since edges $(1,2)$ and $(1,n)$ were chosen without loss of generality as adjacent edges not containing $v$, we get that $\alpha_{e_1}=\alpha_{e_2}$ for any two edges $e_1,e_2$ not containing $v$ by applying this fact to pairs of adjacent edges not containing $v$ repetitively. Let $\alpha_e=A$ for any edge not containing $v$ which we fixed to be $\frac{n}{a-1}$ at the beginning. Then the first equation becomes

$$\left(n-\frac{n}{a-1}\right) \cdot A = \left(n-\frac{n}{a-1} -2\right)\cdot A + \alpha_{(1,\frac{n}{a-1})},$$
which yields that $\alpha_{(1,\frac{n}{a-1})} = 2A$. Since $v=\frac{n}{a-1}$ was chosen without loss of generality, we get that $\alpha_{e'}=2A$ for any edge $e'$ such that $v\in e'$. Note that it is clear that $\alpha_e >0$ for any edge $e$, and thus we can conlude that $\alpha$ is a positive scalar multiple of the left-hand side of the doubling inequality. Since we know the doubling inequality is tight with $T(n,a,2)$ as we've seen through the two types of constructions, it is facet-defining. Finally, note that the sets $S_1$, $S_2$, $S_3$ and $S_4$ all exist only if each part contains at least three vertices, that is, if $n\geq 3(a-1)$.
\qed

\begin{theorem}\label{turandoublinglifting}
The doubling inequality 

$$\sum_{e\in \delta(v)} 2x_e + \sum_{e \in E(K_i)\backslash \delta(v)} x_e \leq ex(i+1,a,2)$$
is facet-defining for any $v\in V(K_i)$ for $T(n,a,2)$ when $i=0 \mod (a-1)$, with $i\geq 3(a-1)$ and $n\geq i$.
\end{theorem}

\proof
By Corollary \ref{lifting}, we only need that for every $e\in E(K_n)\backslash E(K_i)$, there exists an optimal solution $S$ of the doubling inequality in $K_i$ such that $S\cup e$ is Tur\'an. This is clear since, for any such edge $e$, $S\cup e$ is Tur\'an for any optimal solution $S$ of the doubling inequality because any $a$-clique containing the vertices of $e$ as well as some vertices of $K_i$ is missing some edges.
\qed

Clearly, $T(n,a,2)$ has many more facets than the ones we have spoken about. Since the optimal solutions of the Tur\'an graph problem are already known, it's actually quite easy to produce more facets by using proofs like the ones we have seen so far. For example, instead of doubling just one vertex like in the last inequality, we could double two vertices, say $v_1,v_2$. Then 

$$4 \cdot x_{(v_1,v_2)} + 2 \cdot \sum_{\substack{e\in E(K_n):\\ v_1 \text{ or } v_2\in e\\ \text{but not both}}} x_e+ \sum_{\substack{e\in E(K_n):\\ v_1\not\in e, v_2\not\in e}} x_e \leq T(n+2,a,2)$$
is facet-defining if $n=1 \mod (a-1)$, and we can keep playing this game with more vertices. Similarly, it is easy to come up with non-rank facets for the web and wheel graphs. However, these proofs all rely on knowing what optimal solutions look like, knowledge that we are lacking in the next section when considering $r$-uniform hypergraphs with $r\geq 3$. Still, note that if Tur\'an's conjecture for $ex(n,4,3)$ is correct, then our proof for the clique facets for $T(n,a,2)$ could be generalized to $T(n,4,3)$ by using known constructions of extremal graphs for that case (see, e.g., \cite{Frohmader2008}). 

\subsection{Some facets for $T(n,a,r)$}

\subsubsection{Hyperwheel Facets}
Wheel facets for the stable set polytope are related to the graph formed by connecting a single vertex to all vertices of a cycle. We generalize these graphs to hypergraphs and show that they also yield facet-inducing inequalities for $T(n,a,r)$.

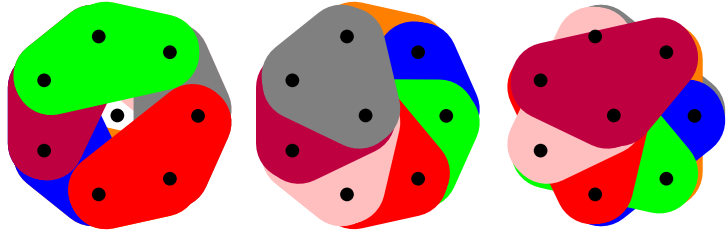
\begin{figure}
\scalebox{0.55}{\begin{tikzpicture}
\node[vertex] (v1) at  (-0.20,4.14) {};
\node[vertex] (v2) at (1.52,3.76) {};
\node[vertex] (v3) at (2.2,2.22) {};
\node[vertex] (v4) at (1.52,0.70) {};
\node[vertex] (v5) at (-0.20,0.32) {};
\node[vertex] (v6) at (-1.52,1.38) {};
\node[vertex] (v7) at (-1.52,3.08) {};
\node [vertex] (v8) at (0.25,2.23) {};

\node[vertex] (v9) at  (5.80,4.14) {};
\node[vertex] (v10) at (7.52,3.76) {};
\node[vertex] (v11) at (8.2,2.22) {};
\node[vertex] (v12) at (7.52,0.70) {};
\node[vertex] (v13) at (5.80,0.32) {};
\node[vertex] (v14) at (4.48,1.38) {};
\node[vertex] (v15) at (4.48,3.08) {};
\node [vertex] (v16) at (6.25,2.23) {};

\node[vertex] (v17) at  (11.80,4.14) {};
\node[vertex] (v18) at (13.52,3.76) {};
\node[vertex] (v19) at (14.2,2.22) {};
\node[vertex] (v20) at (13.52,0.70) {};
\node[vertex] (v21) at (11.80,0.32) {};
\node[vertex] (v22) at (10.48,1.38) {};
\node[vertex] (v23) at (10.48,3.08) {};
\node [vertex] (v24) at (12.25,2.23) {};

\begin{pgfonlayer}{background}
\begin{scope}[transparency group,opacity=0.5]
\draw[edge,opacity=1,color=pink] (v1) -- (v2) -- (v3) -- (v1);
\fill[edge,opacity=1,color=pink] (v1.center) -- (v2.center) -- (v3.center) -- (v1.center);
\end{scope}
\begin{scope}[transparency group,opacity=.5]
\draw[edge,opacity=1,color=gray] (v2) -- (v3) -- (v4) -- (v2);
\fill[edge,opacity=1,color=gray] (v2.center) -- (v3.center) -- (v4.center) -- (v2.center);
\end{scope}
\begin{scope}[transparency group,opacity=.5]
\draw[edge,opacity=1,color=orange] (v5) -- (v6) -- (v4) -- (v5);
\fill[edge,opacity=1,color=orange] (v5.center) -- (v6.center) -- (v4.center) -- (v5.center);
\end{scope}
\begin{scope}[transparency group,opacity=.5]
\draw[edge,opacity=1,color=blue] (v5) -- (v6) -- (v7) -- (v5);
\fill[edge,opacity=1,color=blue] (v5.center) -- (v6.center) -- (v7.center) -- (v5.center);
\end{scope}
\begin{scope}[transparency group,opacity=.5]
\draw[edge,opacity=1,color=purple] (v6) -- (v7) -- (v1) -- (v6);
\fill[edge,opacity=1,color=purple] (v6.center) -- (v7.center) -- (v1.center) -- (v6.center);
\end{scope}
\begin{scope}[transparency group,opacity=.5]
\draw[edge,opacity=1,color=green] (v2) -- (v7) -- (v1) -- (v2);
\fill[edge,opacity=1,color=green] (v2.center) -- (v7.center) -- (v1.center) -- (v2.center);
\end{scope}
\begin{scope}[transparency group,opacity=.8]
\draw[edge,opacity=1,color=red] (v5) -- (v3) -- (v4) -- (v5);
\fill[edge,opacity=1,color=red] (v5.center) -- (v3.center) -- (v4.center) -- (v5.center);
\end{scope}

\begin{scope}[transparency group,opacity=.6]
\draw[edge,opacity=1,color=orange] (v9) -- (v10) -- (v16) -- (v9);
\fill[edge,opacity=1,color=orange] (v9.center) -- (v10.center) -- (v16.center) -- (v9.center);
\end{scope}
\begin{scope}[transparency group,opacity=.55]
\draw[edge,opacity=1,color=blue] (v10) -- (v11) -- (v16) -- (v10);
\fill[edge,opacity=1,color=blue] (v10.center) -- (v11.center) -- (v16.center) -- (v10.center);
\end{scope}
\begin{scope}[transparency group,opacity=.5]
\draw[edge,opacity=1,color=green] (v11) -- (v12) -- (v16) -- (v11);
\fill[edge,opacity=1,color=green] (v11.center) -- (v12.center) -- (v16.center) -- (v11.center);
\end{scope}
\begin{scope}[transparency group,opacity=.8]
\draw[edge,opacity=1,color=red] (v12) -- (v13) -- (v16) -- (v12);
\fill[edge,opacity=1,color=red] (v12.center) -- (v13.center) -- (v16.center) -- (v12.center);
\end{scope}
\begin{scope}[transparency group,opacity=.7]
\draw[edge,opacity=1,color=pink] (v13) -- (v14) -- (v16) -- (v13);
\fill[edge,opacity=1,color=pink] (v13.center) -- (v14.center) -- (v16.center) -- (v13.center);
\end{scope}
\begin{scope}[transparency group,opacity=.7]
\draw[edge,opacity=1,color=purple] (v14) -- (v15) -- (v16) -- (v14);
\fill[edge,opacity=1,color=purple] (v14.center) -- (v15.center) -- (v16.center) -- (v14.center);
\end{scope}
\begin{scope}[transparency group,opacity=.7]
\draw[edge,opacity=1,color=gray] (v9) -- (v15) -- (v16) -- (v9);
\fill[edge,opacity=1,color=gray] (v9.center) -- (v15.center) -- (v16.center) -- (v9.center);
\end{scope}

\begin{scope}[transparency group,opacity=.65]
\draw[edge,opacity=1,color=gray] (v17) -- (v19) -- (v24) -- (v19);
\fill[edge,opacity=1,color=gray] (v17.center) -- (v19.center) -- (v24.center) -- (v19.center);
\end{scope}
\begin{scope}[transparency group,opacity=.6]
\draw[edge,opacity=1,color=orange] (v18) -- (v20) -- (v24) -- (v18);
\fill[edge,opacity=1,color=orange] (v18.center) -- (v20.center) -- (v24.center) -- (v18.center);
\end{scope}
\begin{scope}[transparency group,opacity=.55]
\draw[edge,opacity=1,color=blue] (v19) -- (v21) -- (v24) -- (v19);
\fill[edge,opacity=1,color=blue] (v19.center) -- (v21.center) -- (v24.center) -- (v19.center);
\end{scope}
\begin{scope}[transparency group,opacity=.5]
\draw[edge,opacity=1,color=green] (v20) -- (v22) -- (v24) -- (v20);
\fill[edge,opacity=1,color=green] (v20.center) -- (v22.center) -- (v24.center) -- (v20.center);
\end{scope}
\begin{scope}[transparency group,opacity=.8]
\draw[edge,opacity=1,color=red] (v21) -- (v23) -- (v24) -- (v21);
\fill[edge,opacity=1,color=red] (v21.center) -- (v23.center) -- (v24.center) -- (v21.center);
\end{scope}
\begin{scope}[transparency group,opacity=.7]
\draw[edge,opacity=1,color=pink] (v22) -- (v17) -- (v24) -- (v22);
\fill[edge,opacity=1,color=pink] (v22.center) -- (v17.center) -- (v24.center) -- (v22.center);
\end{scope}
\begin{scope}[transparency group,opacity=.7]
\draw[edge,opacity=1,color=purple] (v23) -- (v18) -- (v24) -- (v23);
\fill[edge,opacity=1,color=purple] (v23.center) -- (v18.center) -- (v24.center) -- (v23.center);
\end{scope}

\end{pgfonlayer}
\end{tikzpicture}}
\caption{A $3$-hyperwheel on $8$ vertices.}
\label{3_W_8_4}
\end{figure}

\begin{definition}
A hyperwheel $^rW_l^a$ is a $r$-uniform hypergraph on $l$ vertices with one vertex in the center, say $l$, and vertices $[l-1]$ placed in a cycle in increasing order around it. The $r$-edges present are such that every $a-1$ consecutive vertices form an $a$-clique with vertex $l$. We only consider wheels for which $n\geq 2a-1$, since otherwise we'd have a complete hypergraph. For example, Figure \ref{3_W_8_4} represents $^3W_8^4$. We call edges that contain the middle vertex, say $l$, \emph{spoke edges} and those who don't, \emph{cycle edges}. Suppose $1 \leq j_1 < \ldots < j_{r-1} \leq l-1$, then we say that the cycle edge $(i, i+j_1 \mod l, \ldots, i+ j_{r-1} \mod l)$ spans $j_{r-1}$ vertices, and the spoke edge $(i, i+j_1 \mod l, \ldots, i+j_{r-2} \mod l, l)$ spans $j_{r-2}$ vertices if $j_{r-1}$ and $j_{r-2}$ are as small as possible.  A spoke edge spans between $r-1$ and $a-1$ vertices of the cycle and a cycle edge, between $r$ and $a-1$ vertices. Note that a $r$-edge spanning $\beta$ vertices of the cycle is contained in $a-1-(\beta-1)=a-\beta$ of the hypercliques. We label the hyperclique spanning vertices $i, i+1, \ldots, i+a-2, l \mod (l-1)$ as \emph{hyperclique $i$}. Note that hypercliques containing any edge $e$ are consecutive modulo $a-1$ in that labeling. The hyperwheel $^rW_l^a$ contains $(l-1)\cdot\binom{a-1}{r-1}$ $r$-edges.
\end{definition}

\begin{theorem}\label{hyperwheelineqs}
The following inequalities are valid and tight for $T(n,a,r)$

$$\sum_{e\in E(^rW_l^a)} x_e \leq |E(^rW_l^a)| - \left\lceil\frac{l-1}{a-r+1}\right\rceil = \binom{a-1}{r-1}\cdot(l-1) - \left\lceil\frac{l-1}{a-r+1}\right\rceil,$$
for all wheels $^rW_l^a$ in $K^r_n$ with $l\leq n$.
\end{theorem}

\proof
We first show that the inequality is valid by determining that it is a Chv\'atal-Gomory cut. Add up the $l-1$ inequalities of the $r$-hypercliques of size $a$ with weight $\frac{1}{a-r+1}$, and the edge inequalities $x_e \leq 1$ for all edges $e$ spanning $\beta$ vertices with weight $\frac{\beta-r+1}{a-r+1}$. This yields the following Chv\'atal-Gomory cut

\begin{multline*}
\sum_{e\in E(^rW_l^a)} x_e\leq \\
\left\lfloor \frac{(l-1)\left(\binom{a}{r} - 1 + \sum_{\beta=r-1}^{a-1} \binom{\beta-2}{\beta-r+1}\cdot (\beta-r+1) + \sum_{\beta=r}^{a-1} \binom{\beta-2}{\beta-r}\cdot(\beta-r+1) \right)}{a-r+1}\right\rfloor
\end{multline*}
since there are $\binom{\beta-2}{\beta-r+1}$ spoke edges starting at any vertex of the cycle and spanning the next $\beta-1$ vertices (where $r-1 \leq \beta \leq a-1$) and there are $\binom{\beta-2}{\beta-r}$ cycle edges starting at any vertex of the cycle and spanning the next $\beta-1$ vertices (where $r \leq \beta \leq a-1$).  We can merge the two sums in the right-hand side together  by recalling that $\binom{\beta-2}{\beta-r+1} + \binom{\beta-2}{\beta-r}= \binom{\beta-1}{\beta-r+1}$ to obtain

$$\left\lfloor \frac{(l-1)\left(\binom{a}{r} - 1 + \sum_{\beta=r}^{a-1} \binom{\beta-1}{\beta-r+1}\cdot (\beta-r+1) \right)}{a-r+1}\right\rfloor$$
because $\binom{r-1-2}{r-1-r+1}\cdot(r-1-r+1)=0$. We now switch the indices to be from $1$ to $a-r$ to simplify the formula:

$$\left\lfloor \frac{(l-1)\left(\binom{a}{r} - 1 + \sum_{\alpha=1}^{a-r} \binom{\alpha+r-2}{\alpha}\cdot \alpha \right)}{a-r+1}\right\rfloor.$$
Then one can easily check that this is equal to $\binom{a-1}{r-1}\cdot(l-1) - \left\lceil\frac{l-1}{a-r+1}\right\rceil$ as desired.

We now show that this inequality is tight by producing an edge set of size $\binom{a-1}{r-1}\cdot(l-1) - \lceil\frac{l-1}{a-r+1}\rceil$ in $^rW_l^a\subset K^r_n$ which contains no hyperclique of size $a$. We want to take all of the edges except a minimum-size set of spoke edges which will ensure that no clique is full. To do so, we remove spoke edges that are contained in as many cliques as possible, namely spoke edges that span only $r-1$ vertices. So first remove such a spoke edge, say $(a-r+1, \ldots, a-1, l)$ without loss of generality, which ensures that the hypercliques starting on vertex $1$ through $a-r+1$ are not full. Then remove spoke edge $(2a-2r+2, \ldots, 2a-r, l)$ which ensures that the hypercliques starting on vertices $a-r+2$ through $2a-2r+2$ are not full, and so on. When the full cycle has been explored that way, we only need to make sure we remove an edge before vertex $1$. Thus, by removing $\left\lceil\frac{l-1}{a-r+1}\right\rceil$ edges, what remains is $a$-hyperclique-free since at least one spoke is missing in each hyperclique. Given that the hyperwheel contains $(l-1)\cdot \binom{a-1}{r-1}$ edges in the first place, it means that such a solution contains 

$$\binom{a-1}{r-1}\cdot(l-1) - \left\lceil\frac{l-1}{a-r+1}\right\rceil$$
$r$-hyperedges and so hyperwheel inequalities of type $^rW_l^a$ are tight for $T(n,a,r)$ with $n \geq l$.
\qed

Note that in the previous proof, we have seen one type of optimal Tur\'an edge set for the hyperwheel $^rW^a_l$ which we call a type I construction. Another type of optimal Tur\'an construction is to remove, without loss of generality, $\lfloor \frac{l-1}{a-r+1} \rfloor$ (r-1)-spanning spoke edges, say edges $(i\cdot (a-r+1),\ldots,i\cdot(a-r+1) + r-2,l)$ for $1\leq i \leq \lfloor \frac{l-1}{a-r+1} \rfloor$ which guarantees that the cliques starting on vertices $1$ through $\lfloor \frac{l-1}{a-r+1} \rfloor$ and spanning the next $a-2$ vertices and central vertex $l$ are all not full. Thus removing any edge contained in all cliques starting on $\lfloor \frac{l-1}{a-r+1} \rfloor$ through $l-1$ will yield an optimal Tur\'an solution which we call of type II.  

\begin{theorem}\label{turanhyperwheelfacet}
The inequality

$$\sum_{e\in E(^rW_l^a)} x_e \leq \binom{a-1}{r-1}\cdot(l-1) - \left\lceil\frac{l-1}{a-r+1}\right\rceil$$
is facet-defining for $T(^rW_l^a, a,2)$ if $l-1 = 1 \mod (a-r+1)$.
\end{theorem}

\proof
Let $\alpha x\leq \beta$ be satisfied by all $x$ in the Tur\'an polytope with $\sum_{e\in E(^rW_l^a)} x_e = \binom{a-1}{r-1}\cdot(l-1) - \left\lceil\frac{l-1}{a-r+1}\right\rceil$; then $\alpha(S)=\beta$ for each Tur\'an edge set $S$ with $|S\cap ^rW_l^a| = \binom{a-1}{r-1}\cdot(l-1) - \left\lceil\frac{l-1}{a-r+1}\right\rceil$.

Consider two distinct spoke edges, including one that spans $r-1$ vertices of the cycle, such that both start from the same vertex, without loss of generality, say $(1,2,\ldots,r-2, l-1,l)$ and $(i_1,\ldots,i_{r-2},l-1,l)$ where $1\leq i_j \leq a-2$ for all $1\leq j \leq r-2$. Let solution $S_1$ be the edge set

\begin{footnotesize}
\begin{multline*}
E(^rW_l^a)\backslash \Bigg\{\left\{(i\cdot (a-r+1),i\cdot (a-r+1)+1,\ldots,i\cdot(a-r+1)+r-2,l)| 1 \leq i \leq \left\lfloor \frac{l-1}{a-r+1} \right\rfloor \right\}\\
\cup (1,2,\ldots,r-2, l-1,l) \Bigg\}
\end{multline*}
\end{footnotesize}
and $S_2$ be

\begin{footnotesize}
\begin{multline*}
E(^rW^a_l)\backslash \Bigg\{\left\{(i\cdot (a-r+1),i\cdot (a-r+1)+1,\ldots,i\cdot(a-r+1)+r-2,l)| 1 \leq i \leq \left\lfloor \frac{l-1}{a-r+1} \right\rfloor \right\} \\
\cup (i_1,\ldots,i_{r-2},l-1,l) \Bigg\},
\end{multline*}
\end{footnotesize}
which are clearly two optimal Tur\'an sets respectively of type I and II in $^rW_l^a$ since $l-1=1 \mod (a-r+1)$ and so the only clique that is still full after removing $\{(i\cdot (a-r+1),i\cdot (a-r+1)+1,\ldots,i\cdot(a-r+1)+r-2,l)| 1 \leq i \leq \lfloor \frac{l-1}{a-r+1} \rfloor \}$ is the clique $(l-1,l,1,2,\ldots,a-2)$ so removing any of those two spoke edges will make the graph $a$-clique-free. This means that $\alpha(S_1)=\beta=\alpha(S_2)$ which implies that $\alpha_{(1,2,\ldots,r-2, l-1,l)}=\alpha_{(i_1,\ldots,i_{r-2},l-1,l)}$. Since we can show this for any spoke edge starting on the same vertex as a spoke edge that spans $r-1$ vertices, we obtain that $\alpha_{e_1}=\alpha_{e_2}$ for any two spoke edges $e_1$ and $e_2$.

Now consider a spoke edge spanning $r-1$ vertices of the cycle and any cycle edge that starts on the same vertex, without loss of generality $(1,2,\ldots,r-2, l-1,l)$ and $(i_1,\ldots,i_{r-1},l-1)$ where $1\leq i_j \leq a-2$ for all $1\leq j \leq r-1$. Again, if we let $S_3$ be the same edge set as $S_1$ and $S_4$ be

\begin{footnotesize}
\begin{multline*}
E(^rW^a_l) \backslash \Bigg\{\left\{(i\cdot (a-r+1),i\cdot (a-r+1)+1,\ldots,i\cdot(a-r+1)+r-2,l)| 1 \leq i \leq \left\lfloor \frac{l-1}{a-r+1} \right\rfloor \right\} \\
\cup (i_1,\ldots,i_{r-1},l) \Bigg\},
\end{multline*}
\end{footnotesize}
then they are both optimal Tur\'an solutions respectively of type I and II by the same argument as we've just seen. Thus, we have that $\alpha(S_3)=\beta=\alpha(S_4)$ which implies that $\alpha_{(1,2,\ldots,r-2, l-1,l)} =\alpha_{(i_1,\ldots,i_{r-1},l-1)}$. Since this is true for any spoke spanning $r-1$ vertices and any cycle edge starting on the same vertex, we have that $\alpha_{e_1}=\alpha_{e_2}$ for any two edges in $^rW_l^a$. It is also clear that $\alpha_e > 0$ for all edges $e\in E(^rW_l^a)$. We thus conclude that $\alpha_3$ is a positive scalar multiple of the left-hand side of the hyperwheel inequality, which is thus facet-defining when $l-1 = 1 \mod (a-r+1)$.
\qed

\begin{theorem}\label{turanhyperwheellifting}
The inequalities

$$\sum_{e\in E(^rW_l^a)} x_e \leq \binom{a-1}{r-1}\cdot(l-1) - \left\lceil\frac{l-1}{a-r+1}\right\rceil$$
are facet-defining for $T(n,a,r)$ for all $^rW_l^a \subseteq K^r_n$ with $l-1 = 1 \mod a-r+1$.
\end{theorem}

\proof
By Theorem \ref{turanliftinghyper01}, we simply need to show that there exists a Tur\'an edge set of size $\binom{a-1}{r-1}\cdot(l-1) - \left\lceil\frac{l-1}{a-r+1}\right\rceil + 1$ in $E(^rW_l^a) \cup e$ for every $e\in E(K^r_n)\backslash E(^rW_l^a)$. First consider an edge $e\in E\backslash E(^rW_l^a)$ that is not an edge of the hyperwheel. If $e=(i_1,\ldots,i_r)$ such that $i_j\not\in V(^rW_l^a)$ for some $1\leq j\leq r$, then it is clear that we can add this edge to any optimal Tur\'an edge set in $^rW_l^a$ without creating a clique of size $a$. Now suppose $e=(i_1,\ldots,i_r)$ such that $i_j \in V(^rW_l^a)$ for all $1\leq j\leq r$. We know that $e$ spans at least $a$ vertices of the cycle since these are the only edges missing. If $e$ spans more than $a$ vertices, then we can add it to any optimal Tur\'an edge set in $^rW_l^a$ without creating an $a$ clique since any $a$-clique containing $e$ would have to contain also another edge spanning at least $a$ vertices, which we know are absent from such an edge set given that they aren't in the hyperwheel in the first place. So we just have to show that there exists an optimal Tur\'an edge set $S$ in $^rW_l^a$ such that $S\cup e$ is still Tur\'an when $e$ spans $a$ vertices.  If $r>2$, then there exist other edges missing in those $a$ vertices, and so they cannot form a clique . When $r=2$, there exists an optimal solution where an edge $(b,b+1)$ is missing, namely the type II construction we discussed with, without loss of generality, let $e=(a-1,l-1)$ and remove  edge $(1,l-1)$ and edges $((a-1)\cdot i,l)$ for $1\leq l \leq \lfloor \frac{l-1}{a-1}\rfloor$ from the wheel. Adding $e$ to this optimal type II Tur\'an solution does not create an $a$-clique.

Thus, if $e\in E(K_n^r)\backslash E(^rW_l^a)$, then there always exists an optimal Tur\'an edge set $S$ in $^rW_l^a$ such that $S\cup e$ is also Tur\'an. Therefore,hyperwheel inequalities on $l$ vertices with $l-1=1 \mod (a-r+1)$ will still be facet-defining for $T(n,a,r)$, $n \geq l$. Actually, by this argument, these wheel inequalities will be facet-inducing for any $T(G,a,r)$ for any graph $G$ that contains such hyperwheels as subgraphs.

\qed

Note that one can give a full linear descriptions of $T(^rW_l^{r+1},r+1,r)$: one needs only clique, non-negativity, edge and hyperwheel inequalities. For a proof, see {\tt http://www.math.washington.edu/$\sim$raymonda/wheel.pdf}. 

\subsubsection{Hyperweb Facets}
A web (or circulant) is a graph with vertices $[n]$ where $(i,j)\in E$ if $i$ and $j$ differ by at most $k$ (mod $n$) and $i\neq j$. Inequalities built from webs are also facet-defining for the stable set polytope. We again consider a generalization to hypergraphs, and show that corresponding inequalities are facet-defining for $T(n,a,r)$.

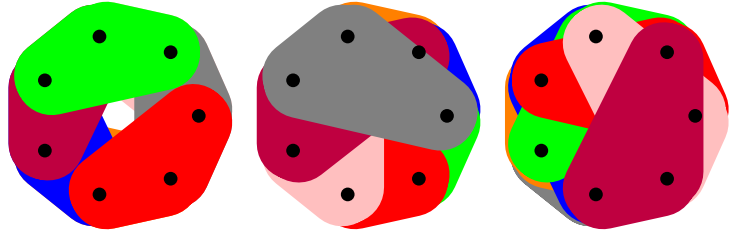
\begin{figure}
\scalebox{0.55}{\begin{tikzpicture}
\node[vertex] (v1) at  (-0.20,4.14) {};
\node[vertex] (v2) at (1.52,3.76) {};
\node[vertex] (v3) at (2.2,2.22) {};
\node[vertex] (v4) at (1.52,0.70) {};
\node[vertex] (v5) at (-0.20,0.32) {};
\node[vertex] (v6) at (-1.52,1.38) {};
\node[vertex] (v7) at (-1.52,3.08) {};

\node[vertex] (v8) at  (5.80,4.14) {};
\node[vertex] (v9) at (7.52,3.76) {};
\node[vertex] (v10) at (8.2,2.22) {};
\node[vertex] (v11) at (7.52,0.70) {};
\node[vertex] (v12) at (5.80,0.32) {};
\node[vertex] (v13) at (4.48,1.38) {};
\node[vertex] (v14) at (4.48,3.08) {};

\node[vertex] (v15) at  (11.80,4.14) {};
\node[vertex] (v16) at (13.52,3.76) {};
\node[vertex] (v17) at (14.2,2.22) {};
\node[vertex] (v18) at (13.52,0.70) {};
\node[vertex] (v19) at (11.80,0.32) {};
\node[vertex] (v20) at (10.48,1.38) {};
\node[vertex] (v21) at (10.48,3.08) {};

\begin{pgfonlayer}{background}
\begin{scope}[transparency group,opacity=0.5]
\draw[edge,opacity=1,color=pink] (v1) -- (v2) -- (v3) -- (v1);
\fill[edge,opacity=1,color=pink] (v1.center) -- (v2.center) -- (v3.center) -- (v1.center);
\end{scope}
\begin{scope}[transparency group,opacity=.5]
\draw[edge,opacity=1,color=gray] (v2) -- (v3) -- (v4) -- (v2);
\fill[edge,opacity=1,color=gray] (v2.center) -- (v3.center) -- (v4.center) -- (v2.center);
\end{scope}
\begin{scope}[transparency group,opacity=.5]
\draw[edge,opacity=1,color=orange] (v5) -- (v6) -- (v4) -- (v5);
\fill[edge,opacity=1,color=orange] (v5.center) -- (v6.center) -- (v4.center) -- (v5.center);
\end{scope}
\begin{scope}[transparency group,opacity=.5]
\draw[edge,opacity=1,color=blue] (v5) -- (v6) -- (v7) -- (v5);
\fill[edge,opacity=1,color=blue] (v5.center) -- (v6.center) -- (v7.center) -- (v5.center);
\end{scope}
\begin{scope}[transparency group,opacity=.5]
\draw[edge,opacity=1,color=purple] (v6) -- (v7) -- (v1) -- (v6);
\fill[edge,opacity=1,color=purple] (v6.center) -- (v7.center) -- (v1.center) -- (v6.center);
\end{scope}
\begin{scope}[transparency group,opacity=.5]
\draw[edge,opacity=1,color=green] (v2) -- (v7) -- (v1) -- (v2);
\fill[edge,opacity=1,color=green] (v2.center) -- (v7.center) -- (v1.center) -- (v2.center);
\end{scope}
\begin{scope}[transparency group,opacity=.8]
\draw[edge,opacity=1,color=red] (v5) -- (v3) -- (v4) -- (v5);
\fill[edge,opacity=1,color=red] (v5.center) -- (v3.center) -- (v4.center) -- (v5.center);
\end{scope}

\begin{scope}[transparency group,opacity=.6]
\draw[edge,opacity=1,color=orange] (v8) -- (v9) -- (v11) -- (v8);
\fill[edge,opacity=1,color=orange] (v8.center) -- (v9.center) -- (v11.center) -- (v8.center);
\end{scope}
\begin{scope}[transparency group,opacity=.55]
\draw[edge,opacity=1,color=blue] (v9) -- (v10) -- (v12) -- (v9);
\fill[edge,opacity=1,color=blue] (v9.center) -- (v10.center) -- (v12.center) -- (v9.center);
\end{scope}
\begin{scope}[transparency group,opacity=.5]
\draw[edge,opacity=1,color=green] (v10) -- (v11) -- (v13) -- (v10);
\fill[edge,opacity=1,color=green] (v10.center) -- (v11.center) -- (v13.center) -- (v10.center);
\end{scope}
\begin{scope}[transparency group,opacity=.8]
\draw[edge,opacity=1,color=red] (v11) -- (v12) -- (v14) -- (v11);
\fill[edge,opacity=1,color=red] (v11.center) -- (v12.center) -- (v14.center) -- (v11.center);
\end{scope}
\begin{scope}[transparency group,opacity=.7]
\draw[edge,opacity=1,color=pink] (v12) -- (v13) -- (v8) -- (v12);
\fill[edge,opacity=1,color=pink] (v12.center) -- (v13.center) -- (v8.center) -- (v12.center);
\end{scope}
\begin{scope}[transparency group,opacity=.7]
\draw[edge,opacity=1,color=purple] (v13) -- (v14) -- (v9) -- (v13);
\fill[edge,opacity=1,color=purple] (v13.center) -- (v14.center) -- (v9.center) -- (v13.center);
\end{scope}
\begin{scope}[transparency group,opacity=.7]
\draw[edge,opacity=1,color=gray] (v14) -- (v8) -- (v10) -- (v14);
\fill[edge,opacity=1,color=gray] (v14.center) -- (v8.center) -- (v10.center) -- (v14.center);
\end{scope}

\begin{scope}[transparency group,opacity=.65]
\draw[edge,opacity=1,color=gray] (v17) -- (v19) -- (v20) -- (v19);
\fill[edge,opacity=1,color=gray] (v17.center) -- (v19.center) -- (v20.center) -- (v19.center);
\end{scope}
\begin{scope}[transparency group,opacity=.6]
\draw[edge,opacity=1,color=orange] (v18) -- (v20) -- (v21) -- (v18);
\fill[edge,opacity=1,color=orange] (v18.center) -- (v20.center) -- (v21.center) -- (v18.center);
\end{scope}
\begin{scope}[transparency group,opacity=.55]
\draw[edge,opacity=1,color=blue] (v19) -- (v21) -- (v15) -- (v19);
\fill[edge,opacity=1,color=blue] (v19.center) -- (v21.center) -- (v15.center) -- (v19.center);
\end{scope}
\begin{scope}[transparency group,opacity=.5]
\draw[edge,opacity=1,color=green] (v20) -- (v15) -- (v16) -- (v20);
\fill[edge,opacity=1,color=green] (v20.center) -- (v15.center) -- (v16.center) -- (v20.center);
\end{scope}
\begin{scope}[transparency group,opacity=.8]
\draw[edge,opacity=1,color=red] (v21) -- (v16) -- (v17) -- (v21);
\fill[edge,opacity=1,color=red] (v21.center) -- (v16.center) -- (v17.center) -- (v21.center);
\end{scope}
\begin{scope}[transparency group,opacity=.7]
\draw[edge,opacity=1,color=pink] (v15) -- (v17) -- (v18) -- (v15);
\fill[edge,opacity=1,color=pink] (v15.center) -- (v17.center) -- (v18.center) -- (v15.center);
\end{scope}
\begin{scope}[transparency group,opacity=.7]
\draw[edge,opacity=1,color=purple] (v16) -- (v18) -- (v19) -- (v16);
\fill[edge,opacity=1,color=purple] (v16.center) -- (v18.center) -- (v19.center) -- (v16.center);
\end{scope}

\end{pgfonlayer}
\end{tikzpicture}}
\caption{A $3$-hyperweb $^3\overline{W}_7^3$ on $7$ vertices.}
\label{3_barW_7_4}
\end{figure}

\begin{definition}
A \emph{hyperweb} $^r\overline{W}_{l}^{a-1}$ is a $r$-uniform hypergraph on $l$ vertices placed in a cycle in increasing order, say $1$ through $l$. The $r$-edges present are such that every $a$ consecutive vertices form an $a$-clique, i.e., for any given vertex, any edge starting with it and spanning at most the next $a-1$ vertices will be present. We only consider hyperwebs for which $l\geq 2a$, so that we do not have a complete hypergraph. A hyperweb is thus for us a hyperwheel with the central vertex removed. For example, Figure \ref{3_barW_7_4} represents the edges present in $^3\overline{W}_7^3$.
\end{definition}

We can prove for hyperwebs theorems similar to the ones we had for hyperwheels. Their proofs are along the same lines, so we just write their statements; the proofs can be found on {\tt http://www.math.washington.edu/$\sim$raymonda/web.pdf}. 

\begin{theorem}
The following inequalities are valid and tight for $T(n,a,r)$

$$\sum_{e\in E(^r\overline{W}_{l}^{a-1})} x_e \leq \binom{a-1}{r-1}\cdot l - \left\lceil\frac{l}{a-r+1}\right\rceil,$$
for all hyperwebs $^r\overline{W}_{l}^{a-1}$ in $K^r_n$ with $l\leq n$.
\end{theorem}

\begin{theorem}\label{turanhyperwebfacet}
Inequality

$$\sum_{e\in E(^r\overline{W}_{l}^{a-1})} x_e \leq \binom{a-1}{r-1}\cdot l - \left\lceil\frac{l}{a-r+1}\right\rceil$$
is facet-defining for $T(^r\overline{W}_{l}^{a-1}, a,2)$ if $l = 1 \mod (a-r+1)$.
\end{theorem}

\begin{theorem}\label{turanhyperweblifting}
The inequalities

$$\sum_{e\in E(^r\overline{W}_{l}^{a-1})} x_e \leq \binom{a-1}{r-1}\cdot l - \left\lceil\frac{l}{a-r+1}\right\rceil$$
are facet-defining for $T(n,a,r)$ for all $^r\overline{W}_{l}^{a-1} \subseteq K^r_n$ with $l=1 \mod (a-r+1)$ and $l \leq n$.
\end{theorem}

\begin{theorem}
We have that

\begin{align*}                                                                                  T(^r\overline{W}_{l}^{r},r+1,r)&=\{x\in \mathbb{R}^{|E(^r\overline{W}_{l}^{r})|} | && x(Q^{r+1}) \leq r & \forall Q^{r+1}\in \mathscr{Q}^{r+1}_{^rW_l^r}\\
&&& x(^r\overline{W}_{l}^{r}) \leq r\cdot l - \left\lceil\frac{l}{2}\right\rceil & \\
&&& 0\leq x_e \leq 1 & \forall e\in E(^r\overline{W}_{l}^{r}) \}\\
\end{align*}
and that all of these inequalities are necessary if $l$ is odd.
\end{theorem}

\section{Conclusion}

Obviously, there is still a lot of work to do on the Tur\'an hypergraph problem and also on its polytope. Hopefully, it is clear that the Tur\'an polytope is interesting in itself, notwithstanding its connection to the famous problem. Its structure is combinatorially rich, and the many parallels that can be drawn between its facets and those of the stable set polytope---one of the polytopes that has been studied the most---make the Tur\'an polytope even more intriguing. In the future, we would love to see whether the rank facets of the stable set polytope for quasi-line graphs (see \cite{Oriolo}) can be transferred to the Tur\'an polytope. 

Understanding some of the facet structure of the Tur\'an polytope also allowed us to understand better why the Tur\'an problem is so hard in general. Indeed, given that the facets we found do not get dominated as the number of vertices grows, this leads us to believe that the number of facets of $T(n,a,r)$ becomes unwieldy as $n$ grows. For example, we proved that cliques of size $i$ were facet-defining for $T(n,3,2)$ for all $i$ odd and $n\geq i$. This means that $T(n,3,2)$ already has at least

$$\binom{n}{3} + \binom{n}{5} + \binom{n}{7} + \ldots + \binom{n}{n'}\approx 2^{n-1}$$
facets, where $n'=n$ if $n$ is odd or $n'=n-1$ if $n$ is even. And these are just the clique facets! Optimizing over a polytope with many facets in a random direction is generally hard, and so this might explain why the problem remains open in general. In another way, we are not trying to optimize in a random direction, so coming up with a clever sequence of Chv\'atal-Gomory cuts might be possible just like it was for the graph case. 

\bibliography{turan}
\bibliographystyle{alpha}
\end{document}